\documentclass[preprint,12pt]{elsarticle}

\usepackage{graphicx} 
\usepackage[ruled,vlined]{algorithm2e}
\usepackage{bm}
\usepackage{algorithmic}
\usepackage{amssymb} 
\usepackage{amsmath}
\usepackage{amsthm}
\usepackage{url}
\usepackage{caption}
\usepackage{subcaption}
\usepackage{multirow}

\usepackage{xcolor} 

\begin{document}

\begin{frontmatter}

\title{Enhanced fifth order WENO Shock-Capturing Schemes \\ with Deep Learning}
\author{Tatiana Kossaczká\corref{cor1}}
\ead{kossaczka@uni-wuppertal.de}
\author{Matthias Ehrhardt} 
\ead{ehrhardt@uni-wuppertal.de}
\author{Michael Günther}
\ead{guenther@uni-wuppertal.de}
\address{Institute of Mathematical Modelling, Analysis and Computational Mathematics (IMACM), 
Chair of Applied Mathematics and Numerical Analysis, Bergische Universität Wuppertal, Gaußstraße 20, 42119 Wuppertal, Germany}
\cortext[cor1]{corresponding author}

\begin{abstract}
In this paper we enhance the well-known fifth order WENO shock-capturing scheme 
by using deep learning techniques.
This fine-tuning of an existing algorithm is implemented by training a rather small neural network to modify the smoothness indicators of the WENO scheme in order to improve the numerical results especially at discontinuities. 
In our approach no further post-processing is needed to ensure
the consistency of the method, which simplifies the method and increases the effect of the neural network. Moreover, the convergence of the resulting scheme can be theoretically proven.
 
We demonstrate our findings with the inviscid Burgers’ equation, the Buckley-Leverett equation and the 1-D Euler equations of gas dynamics.
Hereby we investigate the classical Sod problem and the Lax problem and show that our novel method outperforms the classical fifth order WENO schemes in simulations where the numerical solution is too diffusive or tends to overshoot at shocks.
\end{abstract}


\end{frontmatter}


\section{Introduction}\label{S:1}
Typically, numerical fluid mechanics deals with nonlinear hyperbolic partial differential equations (PDEs). 
In its simplest one-dimensional form, these equations can be represented as
\begin{equation} \label{eq:hyp}
     \frac{\partial u}{\partial t}+\frac{\partial f(u)}{\partial x} = 0,\qquad t>0,
\end{equation}
where $x$ represents space, $t$ denotes time, $u(x,t)$ is conserved quantity and $f(u(x,t))$ is its flux.
Discontinuous initial data develops after finite time 
a discontinuity ('shock wave') or a rarefaction wave
(regardless of the smoothness of the initial or boundary data). 
Hence, suitable numerical methods for these problems must be designed to properly handle especially the accurate solution of the discontinuities.
First in 1980, Crandall and Majda \cite{Crandall} proposed the class of \textit{monotone schemes}
that are nonlinearly stable in the $L_1$ norm
and satisfy certain entropy conditions. It can be proven that the corresponding solutions converge
to bounded variation entropy solutions including error estimates.
However, these schemes are only first order accurate
and by the fundamental Godunov theorem \cite{Godunov} it is known that one has to consider nonlinear non-oscillatory schemes
to overcome this accuracy barrier.

In this direction so-called \textit{shock-capturing schemes} were developed that were able to resolve sharply a shock or a sharp gradient region without introducing too much diffusion or overshoot behaviour \cite{Harten}. 
Additionally, to remedy the above mentioned drawback, at regions with smooth flow these schemes exhibit a rather high order of convergence.
The well-known representative of this class of methods are the \textit{essentially non-oscillatory (ENO) schemes} \cite{Harten87}
with high order accuracy in smooth regions
and sharply resolving shocks in an essentially
non-oscillatory way using a smoothness indicator function,
see e.g.\ \cite{Shu1998}.
Later on Jiang and Shu \cite{Jiang} further improved these schemes and proposed a weighted ENO (in the sequel abbreviated with WENO-JS) scheme,
that is still regarded as a state-of-the-art solution approach. 

Subsequently, different new strategies were developed by modifying the WENO-JS schemes, i.e.\ by altering by smoothness indicators \cite{Borges, Castro, Kim, Rathan, Ha, Li20} or by modifying the nonlinear weights \cite{Henrick}.
Besides, another goal in optimizing these schemes was to
minimize the dispersion error (dispersion-relation-preserving (DRP) schemes) \cite{Liu13, Tam93}, also combined with the WENO
approach leading to optimized WENO (OWENO) schemes \cite{Wang01}.

Recently, machine learning was widely used to compute the solution of PDEs. We refer to \cite{lagaris, Sirignano, Berg}, where the neural network algorithm is used to approximate a solution of a particular PDE problem. 
Following that approach, the solution of a particular PDE is a result of a neural network training procedure.
Another idea is to improve a specific numerical scheme using neural networks. 
The training of a neural network is made offline and results in a new numerical scheme applicable to a wider class of problems. This idea was recently used by Beck et al.\ \cite{Beck} for discontinuous Galerkin methods or in \cite{hsieh2019learning} for learning iterative PDE solvers and we also follow this approach. 

The recent work of Stevens and Colonius \cite{Stevens} introduces new WENO-NN scheme based on neural network algorithm. 
In their work, the finite-volume coefficients of the WENO-JS scheme are perturbed, while maintaining the original smoothness indicators and nonlinear weights. However, the resulting scheme presented in their paper has only first order of convergence. Another neural network based WENO scheme was developed by Liu and Wen \cite{Liu20}, where the new smoothness indicators are an output of the neural network algorithm.

We implement in our work another WENO extension 
based on deep learning. This approach is used to improve the classic WENO-JS and 
WENO-Z \cite{Borges} scheme in this paper, but could be efficiently applied also to other WENO modifications.
For this purpose we will train a rather small neural network to perturb the smoothness indicator functions of the WENO-JS scheme. 
As we do not develop any new smoothness indicators as in \cite{Liu20}, but only their multiplicative perturbations, 
we are able to prove the convergence of the resulting scheme. 
We call this new scheme WENO-DS (Deep Smoothness), as we modify the smoothness indicators using deep neural networks. 
This scheme has less diffusion and less overshoot in shocks than the WENO-JS and the WENO-Z scheme, while maintaining high order accuracy in smooth regions. 
Unlike the recent work of Stevens and Colonius~\cite{Stevens}, we do not need any post-processing steps, 
which simplifies the procedure and also increases the effect of the deep learning part. 

Finally, let us note that the use of WENO methods is not limited to hyperbolic PDEs, see e.g.\
\cite{Tatiana:Thesis:2019} for an application in finance.

The paper is structured as follows. 
In Section~\ref{S:2} we introduce the WENO-JS and WENO-Z schemes under consideration in detail.
Next, in Section~\ref{S:22} we introduce our deep learning approach, where neural networks are used to further improve WENO methods without any post-processing.
The corresponding convergence proofs for two WENO schemes are given in Section~\ref{sec:convergence}.
Then we present our numerical results in Section~\ref{S:num}, which illustrate the improvements of our proposed method.
Finally, we conclude our work in Section~\ref{concl}.

\section{The WENO Scheme}\label{S:2}
Let $\{I_i\}$ be the partition of a spatial domain with the $i$-th cell $I_i=[x_{i-\frac{1}{2}}, x_{i+\frac{1}{2}}]$. 
We consider a uniform grid defined by the points $x_i=x_0+i\Delta x$, $i=0,\dots,N$, which are the centers of the cells with cell boundaries defined by $x_{i+\frac{1}{2}}=x_i+\frac{\Delta x}{2}$. 
The value of a function $f$ at $x_i$ is indicated by $f_i=f(x_i)$.

The spatial discretization of one-dimensional conservation laws \eqref{eq:hyp} yields a system of ordinary differential equations ('method of lines') and the resulting semi-discrete scheme is
\begin{equation} \label{eq:hypsemi}
    \frac{du_i}{dt} = -\frac{1}{\Delta x}(\hat{f}_{i+\frac{1}{2}}-\hat{f}_{i-\frac{1}{2}}),
\end{equation}
where $\hat{f}$ is a numerical approximation of the flux function $f$.
Following \cite{Jiang}, if we define a function $h$ implicitly by 
\begin{equation} \label{eq:int}
    f\bigl(u(x)\bigr) = \frac{1}{\Delta x} \int_{x-\frac{\Delta x}{2}}^{x+\frac{\Delta x}{2}} h( \xi)\, d\xi,
\end{equation}
then \eqref{eq:hypsemi} is approximated by
\begin{equation} \label{eq:approx}
    f'\bigl(u(x_i)\bigr) = \frac{1}{\Delta x}\bigl(h_{i+\frac{1}{2}} - h_{i-\frac{1}{2}}\bigr),
    \qquad h_{i \pm \frac{1}{2}} = h(x_{i \pm \frac{1}{2}}),
\end{equation}
where $h_{i \pm \frac{1}{2}}$ approximates the numerical flux $\hat{f}_{ \pm \frac{1}{2}}$ with fifth order of accuracy
\begin{equation} \label{eq:flux_fifth}
    \hat{f}_{ \pm \frac{1}{2}} = h_{i \pm \frac{1}{2}} + O(\Delta x^5).
\end{equation}
This results in a conservative numerical scheme.

To guarantee the stability of the method and avoid entropy violating solutions, 
the \textit{flux splitting method} is applied, thus
\begin{equation} \label{eq:fluxsplit}
    f(u) = f^+(u) + f^-(u),\quad
    \text{where}\quad \frac{d f^+(u)}{du}\ge0\quad
    \text{and}\quad \frac{d f^-(u)}{du}\le0.
\end{equation}
The numerical flux $\hat{f}_{i \pm \frac{1}{2}}$ is then represented by 
$\hat{f}_{i \pm \frac{1}{2}} = \hat{f}_{i \pm \frac{1}{2}}^+ + \hat{f}_{i \pm \frac{1}{2}}^-$ and the final scheme is formed as
\begin{equation} \label{eq:hypsemi_splitting}
    \frac{du_i}{dt} = -\frac{1}{\Delta x}\left[\left(\hat{f}_{i+\frac{1}{2}}^+-\hat{f}_{i-\frac{1}{2}}^+\right)+\left(\hat{f}_{i+\frac{1}{2}}^--\hat{f}_{i-\frac{1}{2}}^-\right)\right].
\end{equation}
Next, only  the construction of $\hat{f}_{i \pm \frac{1}{2}}^+$ is considered. 
The negative part can be then obtained using symmetry (see e.g.\ \cite{Wang}). 

\subsection{Fifth order WENO scheme}
\label{S:20}
For a construction of $\hat{f}_{i + \frac{1}{2}}$, the fifth order WENO method uses a 5-point stencil
\begin{equation} \label{eq:stencil}
    S(i)=\{x_{i-2}, \dots,x_{i+2} \}
\end{equation}
divided into three candidate substencils, which are 
given by
\begin{equation} \label{eq:substencils}
    S^m(i)=\{x_{i+m-2}, x_{i+m-1}, x_{i+m} \}, \quad  m=0,1,2.
\end{equation}
To form the numerical flux over the entire 5-point stencil, the numerical flux 
for each of these substencils $\hat{f}^m_{i+\frac{1}{2}} = h_{i+\frac{1}{2}}+ O(\Delta x^3)$  is calculated. 
These fluxes are then averaged in such a way, that fifth order convergence is ensured in the smooth regions. 
In regions with discontinuities, the weights should partly remove the contribution of these stencils 
so that the solution near the shock can be approximated in more stable manner.

Let $\hat{f}^m(x)$ be the polynomial approximation of $h(x)$ on each of the substencils \eqref{eq:substencils}. 
Then, evaluated at $i + \frac{1}{2}$ we obtain
\begin{equation} \label{eq:flux_plus}
   \hat{f}^m(x_{i+\frac{1}{2}}) =  \hat{f}^m_{i+\frac{1}{2}} 
   = \sum_{j=0}^{2} c_{m,j}\, f_{i+m-2+j}, \quad i=0,\dots,N,
\end{equation}
where $c_{m,j}$ are the Lagrangian interpolation coefficients, dependent on $m$ (see \cite{Jiang}).
They take an explicit form
\begin{align}
\hat{f}^0_{i + \frac{1}{2}} &=  \frac{2f(u_{i-2})-7f(u_{i-1})+11f(u_{i})}{6}, \nonumber \\
\hat{f}^1_{i + \frac{1}{2}} &=  \frac{-f(u_{i-1})+5f(u_{i})+2f(u_{i+1})}{6}, \\
\hat{f}^2_{i + \frac{1}{2}} &=  \frac{2f(u_{i})+5f(u_{i+1})-f(u_{i+2})}{6}  \nonumber
\end{align}
and the numerical fluxes $\hat{f}^m_{i - \frac{1}{2}}$ can be obtained by shifting each index by $-1$.
Using the Taylor series expansion it can be shown that:
\begin{equation} 
     \hat{f}^m_{i\pm\frac{1}{2}} = h_{i\pm\frac{1}{2}} +A_m \Delta x^3 
     + O(\Delta x^4).
\end{equation} 

Then, the convex combination of the interpolated values $\hat{f}^m(x_{i \pm \frac{1}{2}})$ given by 
\begin{equation}\label{eq:omegam}
   \hat{f}_{i\pm\frac{1}{2}}
   = \sum_{j=0}^{2} \omega_m \,\hat{f}^m(x_{i \pm \frac{1}{2}})
\end{equation}
yields the WENO approximation of the value $h_{i \pm \frac{1}{2}}$, 
where $\omega_m$ are the nonlinear weights defined as, cf.\  \cite{Jiang} 
\begin{equation} \label{eq:omegas}
   \omega_m^{JS} = \frac{\alpha_m^{JS}}{\sum_{i=0}^{2}\alpha_i^{JS}}, 
   \quad \text{ where } \quad  
   \alpha_m^{JS} = \frac{d_m}{ (\epsilon + \beta_m)^2 }.
\end{equation}
The scheme using these nonlinear weights is denoted as the WENO-JS scheme. 
The parameter $\epsilon$ guarantees that the denominator does not become zero, and the coefficients $d_0$, $d_1$ and $d_2$ are called ideal weights, 
which would form the upstream fifth order central scheme for the 5-point stencil and satisfy \eqref{eq:flux_fifth}.
Their values are:
\begin{equation} 
    d_0=\frac{1}{10}, \quad d_1=\frac{6}{10}, \quad d_2=\frac{3}{10}.
\end{equation}
The parameter $\beta_m$ is called the \textit{smoothness indicator} and is analyzed in the next section.

\subsection{Smoothness Indicators} 
\label{S:21}
The role of smoothness indicators is to measure the regularity of the polynomial approximation  $\hat{f}^m(x)$ in each of three substencils.
As developed in \cite{Jiang}, they are defined as: 
\begin{equation} \label{eq:beta}
   \beta_m = \sum_{q=1}^{2} \Delta x^{2q-1} \int_{x_{i-\frac{1}{2}}}^{x_{i+\frac{1}{2}}} 
      \Bigl(\frac{d^q \hat{f}^m(x)}{dx^q} \Bigr)^2\,dx.
\end{equation} 
Corresponding to the flux approximation $\hat{f}_{i+\frac{1}{2}}$ they take an explicit form
\begin{align}
\beta_0 &= \frac{13}{12} \bigl(f(u_{i-2})-2f(u_{i-1})+f(u_{i})\bigr)^2 +\frac{1}{4}\bigl(f(u_{i-2})-4f(u_{i-1})+3f(u_{i})\bigr)^2 , \nonumber \\
\beta_1 &= \frac{13}{12} \bigl(f(u_{i-1})-2f(u_{i})+f(u_{i+1})\bigr)^2 +\frac{1}{4}\bigl(-f(u_{i-1})+f(u_{i+1})\bigr)^2, \\
\beta_2 &= \frac{13}{12} \bigl(f(u_{i})-2f(u_{i+1})+f(u_{i+2})\bigr)^2 +\frac{1}{4}\bigl(3f(u_{i})-4f(u_{i+1})+f(u_{i+2})\bigr)^2  \nonumber 
\end{align}
and their Taylor expansions at $x_i$ are:
\begin{equation}\label{eq:taylor_beta}
    \begin{split}
\beta_0 &= f_{x}^2 \Delta x^2 + \Bigl( \frac{13}{12} f_{xx}^2 - \frac{2}{3} f_{x}f_{xxx} \Bigr) \Delta x^4 \\
&\qquad + \Bigl(-\frac{13}{6} f_{xx}f_{xxx} + \frac{1}{2} f_{x}f_{xxxx} \Bigr) \Delta x^5 + O(\Delta x^6), \\
\beta_1 &= f_{x}^2 \Delta x^2 + \Bigl( \frac{13}{12} f_{xx}^2 + \frac{1}{3} f_{x}f_{xxx} \Bigr) \Delta x^4 + O(\Delta x^6), \\ 
\beta_2 &= f_{x}^2 \Delta x^2 + \Bigl( \frac{13}{12} f_{xx}^2 - \frac{2}{3} f_{x}f_{xxx} \Bigr) \Delta x^4\\
&\qquad+ \Bigl(\frac{13}{6} f_{xx}f_{xxx} - \frac{1}{2} f_{x}f_{xxxx} \Bigr) \Delta x^5 + O(\Delta x^6). 
\end{split} 
\end{equation} 
These indicators are designed to come closer to zero for smooth parts of the solution so that the nonlinear weights $\omega_m$ come closer to the ideal weights $d_m$. 
In the case that the stencil $S^m$ contains a discontinuity, $\beta_m$ is $O(1)$ and the corresponding weight $\omega_m$ becomes smaller, therefore the contribution of the substencil $S^m$  is reduced. 

Following \cite{Henrick}, it can be shown that demanding \eqref{eq:flux_fifth} we obtain the sufficient conditions for the fifth order convergence:
\begin{align}
  \sum_{m=0}^{2}(\omega_m^{\pm} - d_m) &= O(\Delta x^6), \label{eq:cond1_0} \\
       \omega_m^{\pm} - d_m &= O(\Delta x^3). \label{eq:suff_cond}
\end{align}
Considering the overall finite difference formula $\hat{f}_{j+\frac{1}{2}} - \hat{f}_{j-\frac{1}{2}} = f'(x) \Delta x + O(\Delta x^6)$, it can be shown, that \eqref{eq:suff_cond} may be relaxed and we obtain the following sufficient and necessary conditions:
\begin{align}
  \sum_{m=0}^{2}(\omega_m^{\pm} - d_m) &= O(\Delta x^6), \label{eq:cond1} \\
   \sum_{m=0}^{2}A_m(\omega_m^{+} - \omega_m^{-}) &= O(\Delta x^3),  \label{eq:cond2}  \\
  \omega_m^{\pm} - d_m &= O(\Delta x^2). \label{eq:cond3} 
\end{align}
Note that due to the normalization \eqref{eq:omegas}, the first condition \eqref{eq:cond1} (resp. \eqref{eq:cond1_0}) is always fulfilled. (The superscripts $\pm$ on $\omega_m$ specify their use in $\hat{f}_{i + \frac{1}{2}}$ or $\hat{f}_{i - \frac{1}{2}}$).

The convergence analysis was performed in \cite{Jiang} and it was shown that if 
\begin{equation} \label{eq:betaD}
    \beta_m = D\bigl(1+O(\Delta x^{2})\bigr),
\end{equation}
with $D$ being a non-zero constant independent of $m$, 
the condition \eqref{eq:cond3} is satisfied and the scheme has the expected fifth order accuracy.
However, it was shown in \cite{Henrick} that at the critical points where the first derivative of $f$ vanishes, 
the convergence order of the scheme from \cite{Jiang} decreases to the third order. 
Moreover, if the second derivative also vanishes, the order is further reduced to the second order. 
For a further explanation of this problem we refer to \cite{Henrick}.

\subsection{The WENO-Z scheme}
In this paper we consider the modified WENO scheme of Borges et al.\ \cite{Borges} with a new global smoothness indicator,
which is characterized by 
\begin{equation} \label{eq:tau}
    \tau_5 = |\beta_0 - \beta_2|.
\end{equation}
It is easy to see from the equations \eqref{eq:taylor_beta} that
\begin{equation}
    \tau_5 = \frac{13}{3}|f_{xx}f_{xxx}|\Delta x^5 + O(\Delta x^6).
\end{equation}
The new WENO-Z weights are then defined by
\begin{equation} \label{eq:WENOZ}
   \omega_m^Z =  \frac{\alpha^{Z}_m}{\sum_{i=0}^{2}\alpha^{Z}_i}, \quad \text{ where } \quad \alpha^{Z}_m
    = d_m \left[ 1+ \left( \frac{\tau_5}{\beta_m + \epsilon} \right)^2 \right].
\end{equation}
Borges et al.\ \cite{Borges} have shown that when using these nonlinear weights, fifth order convergence is preserved, even at the critical points where $f'(u)=0$.

\section{The Deep Learning Approach for WENO Schemes}\label{S:22}
To better capture discontinuities and avoid oscillations, we propose to apply deep learning to develop new smoothness indicators. 
We construct them as products of the original smoothness indicators $\beta_m$ and multipliers $\delta_m$ which are outputs of a neural network algorithm. 
We refer to these new smoothness indicators as $\beta_m^{DS}$, where index $DS$ corresponds to the new WENO-DS scheme:
\begin{equation} \label{eq:betann}
    \beta_m^{DS} = \beta_m (\delta_m + C),
\end{equation}
where $C$ is a constant, which role is crucial for the proof of consistency and convergence and we will explain how to choose it in the Section~\ref{sec:convergence}.
We emphasize that this formulation is very advantageous in a sense that the consistency and convergence 
properties of the original WENO method are preserved. 
In the case that the solution is smooth and the original smoothness indicator $\beta_m$ converges to zero, 
the smoothness indicator $\beta_m^{DS}$ behaves in the same way. 
If the smoothness indicator $\beta_m$ is $O(1)$, the multiplier $\delta_m$ can change it so that the final scheme performs better.

\begin{figure}[ht!] 
\centering\includegraphics[width=1\linewidth]{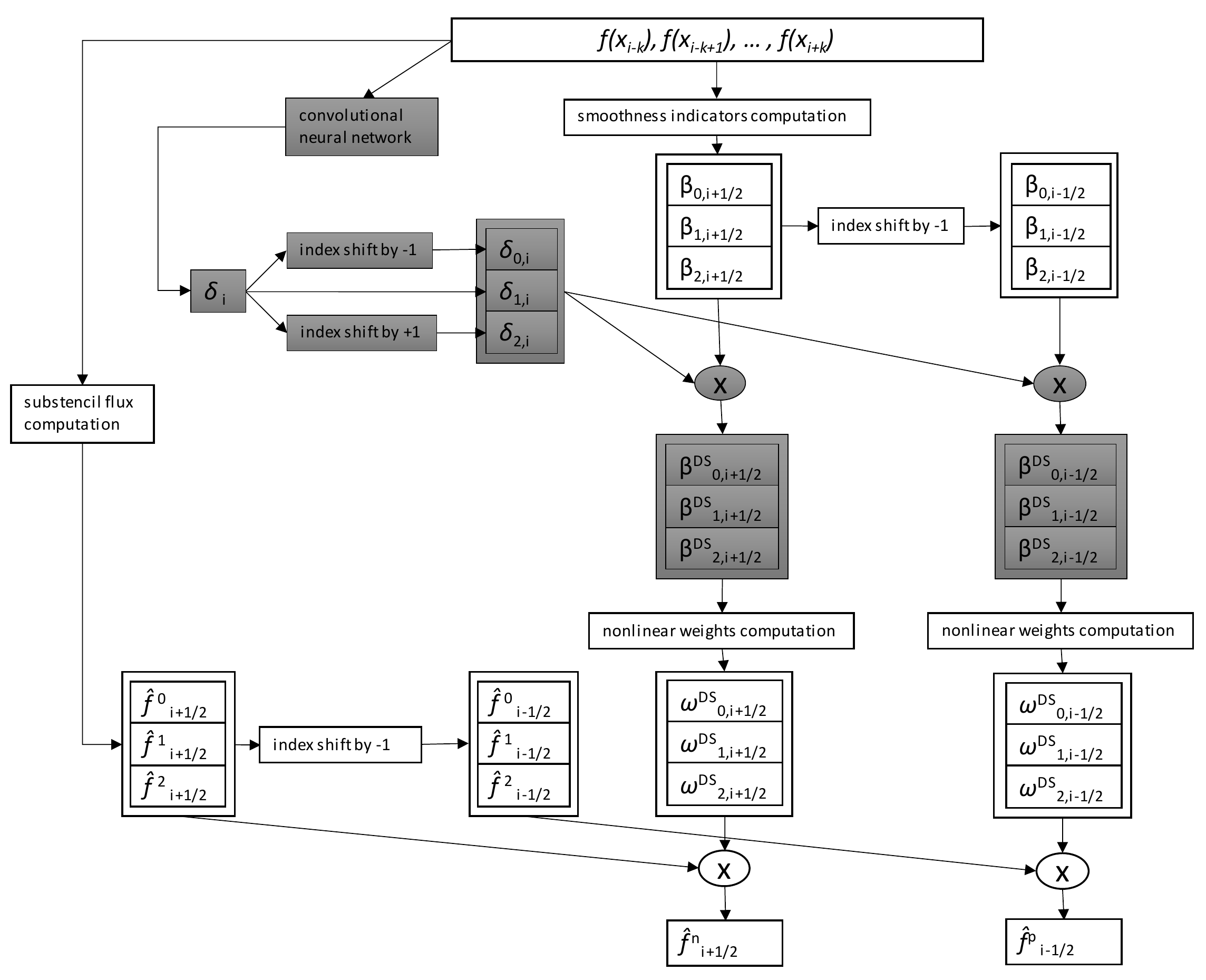}
\caption{The structure of the WENO method combined with the neural network algorithm. 
The white parts of the graph correspond to the original WENO method. 
The grey parts are added to this method so that the whole graph corresponds to the new method WENO-DS. $2k+1$ is the size of the receptive field of the whole CNN, x denotes the element-wise multiplication.}
\label{fig:WENO_whole} 
\end{figure}

In the original WENO method, the stencil \eqref{eq:stencil} is used to approximate the solution in $x_i$, 
and the fluxes are being reconstructed in the points $x_{i-\frac{1}{2}}$ and $x_{i+\frac{1}{2}}$.
To define $\hat{f}^m_{i-\frac{1}{2}}$ we use \eqref{eq:flux_plus} and shift each index by $-1$. 
In our approach we proceed as in the classical WENO method \cite{Shu1998} and compute the smoothness indicators 
as described in \eqref{eq:beta} in the Section~\ref{S:21}.
We use them for a flux reconstruction $\hat{f}_{i+\frac{1}{2}}$ 
and then by shifting each of the index by $-1$ we compute the smoothness indicators corresponding to the 
flux approximation $\hat{f}_{i-\frac{1}{2}}$. 
We denote them as $\beta_{m,i+\frac{1}{2}}$ and $\beta_{m,i-\frac{1}{2}}$, respectively. 
For a fixed $m$ we could make the multiplier $\delta_m$ for $\beta_{m,i+\frac{1}{2}}$ and $\beta_{m,i-\frac{1}{2}}$ dependent on the location of the substencils corresponding to $\beta_{m,i+\frac{1}{2}}$ and $\beta_{m,i-\frac{1}{2}}$. 
This would result in two different multipliers for $\beta_{m,i+\frac{1}{2}}$ and $\beta_{m,i-\frac{1}{2}}$. 
However, experimentally we got the superior results by using the same multiplier $\delta_{m,i}$ for both $\beta_{m,i+\frac{1}{2}}$ and $\beta_{m,i-\frac{1}{2}}$, dependent only on the position $i$ of the global stencil.
The new smoothness indicators are then computed as
\begin{equation}
\begin{aligned} \label{eq:betas_DS_i}
     \beta_{m,i+\frac{1}{2}}^{DS} &= \beta_{m,i+\frac{1}{2}} (\delta_{m,i} + C), \\
          \beta_{m,i-\frac{1}{2}}^{DS} &= \beta_{m,i-\frac{1}{2}} (\delta_{m,i} + C) 
\end{aligned}
\end{equation}
and the values $\delta_0$, $\delta_1$, $\delta_2$ are obtained by simple index shift so that it holds
\begin{equation} \label{eq:delta_shift}
    \delta_{0,i+1} =  \delta_{1,i} =  \delta_{2,i-1}, \quad i=0, \ldots, N.
\end{equation}
Finally, we obtain the flux approximations in the same way as in classical WENO schemes, but using the new smoothness indicators \eqref{eq:betas_DS_i}:
\begin{equation} \label{eq:flux_pn}
    \hat{f}^p_{i-\frac{1}{2}} \quad \text{and} \quad \hat{f}^n_{i+\frac{1}{2}},
\end{equation}
which are used for approximating of a solution in a point $x_i$. Superscripts $p$ and $n$ indicate the difference between the values at the same location $x_{i+\frac{1}{2}}$  when we approximate the solution in points $x_i$ and $x_{i+1}$ (resp.\ at the same location $x_{i-\frac{1}{2}}$  when we approximate the solution in  points $x_{i-1}$ and $x_{i}$). We present the whole algorithm of the method in the Figure~\ref{fig:WENO_whole}. 

As we mentioned before, the flux splitting technique \eqref{eq:fluxsplit} is used. Each part of a flux, $f^+$ and $f^-$ represents different type of input data to the neural network. 
Therefore we use two neural networks, for the positive and negative part of a flux with the input values $f^+(x_i)$ to the first neural network and  $f^-(x_i)$ to the second neural network, $i=0, \ldots, N$.
Each of the neural networks produces different outputs, multipliers corresponding to the positive and negative part of a flux.
For simplicity we further drop the superscripts $\pm$ and when we talk about the inputs to the neural network we always mean both $f^+(x_i)$ and  $f^-(x_i)$.  We denote by $\bar{f}(x)$ the vector $(f(x_0), f(x_1), \ldots f(x_N))$ and formulate the neural network as a function $F\bigl(\bar{f}(x)\bigr)$. 

To ensure consistency we propose the use of a \textit{convolutional neural network (CNN)}. Firstly, this type of neural network is computationally efficient and secondly, 
it makes the multipliers independent of their position in the spatial grid so that the final numerical scheme is spatially invariant. 
To ensure the convergence of the method, we use the differentiable activation functions like the exponential linear unit (ELU) and the sigmoid function. 
If the layers of the neural network are differentiable functions, then their composition, the neural network function $F(\cdot)$, is also a differentiable function. 
The ELU activation function has the advantage that it does not cause the dying gradient problem, 
the sigmoid activation function ensures that the output of the neural network is between 0 and 1.

Let  us note that we use for the implementation Python with the deep learning library PyTorch \cite{NEURIPS2019_9015} (\url{https://pytorch.org/}), which is capable of GPU acceleration.

\section{Convergence analysis of the new WENO scheme (WENO-DS)}
\label{sec:convergence}

\subsection{Convergence analysis of WENO-JS scheme with new smoothness indicators $\beta_{m,i \pm \frac{1}{2}}^{DS}$} \label{S:41}

Let us express the multipliers $\delta_{m,i}$ for the smoothness indicators $\beta_{m,i \pm \frac{1}{2}}$, $m=0,1,2$ used in the node $x_i$ as the outputs of a neural network function.
Following \eqref{eq:delta_shift} and using the fact, that the layers of the CNN are spatially invariant differentiable functions, we can write
\begin{equation}
\begin{split}
  \delta_{0,i} &= F\bigl(\bar{f}(\bar{x}_{i-1})\bigr) = \Phi(\bar{x}_{i}-\Delta x) = \Phi(\bar{x}_{i}) - O(\Delta x), \\
  \delta_{1,i} &=F\bigl(\bar{f}(\bar{x}_{i})\bigr) = \Phi(\bar{x}_{i}), \\
  \delta_{2,i} &=F\bigl(\bar{f}(\bar{x}_{i+1})\bigr) = \Phi(\bar{x}_{i}+\Delta x) = \Phi(\bar{x}_{i}) + O(\Delta x), 
\end{split}\end{equation}
where 
\begin{equation}
\begin{split}
\bar{x}_{i}&=(x_{i-k}, x_{i-k+1}, \ldots, x_{i+k}), \\
\bar{f}(\bar{x}_i) &= (f(x_{i-k}), f(x_{i-k+1}), \ldots, f(x_{i+k})),
\end{split}
\end{equation}
where $2k+1$ is the size of the receptive field of the whole CNN and  $\Phi$ is the function composition $F \circ \bar{f}$. 
Then using \eqref{eq:betaD} it holds
\begin{equation} \label{eq:new_beta}
    \beta_{m,i \pm \frac{1}{2}}^{DS} = \beta_{m,i \pm \frac{1}{2}} (\delta_{m,i} + C) 
    =  D\bigl(1+O(\Delta x^{2})\bigr) \bigl(\Phi(\bar{x}_{i}) + O(\Delta x) + C\bigr).
\end{equation}
We denote $P(\bar{x}_{i})=\Phi(\bar{x}_{i})+C$ and we set $C$ such that $\Phi(\bar{x}_{i})+C>\kappa>0$ with $\kappa$ fixed. 
Then we ensure that $P(\bar{x}_{i})=O(1)$. 
Performing the multiplication in \eqref{eq:new_beta} we obtain
\begin{equation} \label{eq:beta_expand}
\begin{split}
  \beta_{m,i \pm \frac{1}{2}}^{DS} &= D \bigl(P(\bar{x}_{i}) +P(\bar{x}_{i})O(\Delta x^2) + O(\Delta x) + O(\Delta x^3) \bigr)  \\
  &= DP(\bar{x}_{i}) \bigl( 1 + O(\Delta x) \bigr) = \Tilde{D}\bigl(1 + O(\Delta x) \bigr).
\end{split}\end{equation}
Here we can proceed as in \cite{Henrick}, but for the reader's convenience we repeat the steps of the proof: 
insert \eqref{eq:beta_expand} into \eqref{eq:omegas} and take $\epsilon=0$
\begin{equation} 
  \alpha_{m,i \pm \frac{1}{2}}^{DS}= \frac{d_m}{\bigl(\Tilde{D}(1+O(\Delta x)\bigr)^2} = \frac{d_m}{\Tilde{D}^2}\bigl(1+O(\Delta x)\bigr).
\end{equation}
This implies that 
\begin{equation}
   \sum_{m=0}^{2}\alpha_{m,i \pm \frac{1}{2}}^{DS} = \frac{1}{\Tilde{D}^2}\bigl(1+O(\Delta x)\bigr),
\end{equation}
where we used the fact that $\sum_{m=0}^2 d_m = 1$. Finally, substituting into \eqref{eq:omegas} we obtain
\begin{equation} \label{eq:WENO_DS_JS_cond}
   \omega_{m,i \pm \frac{1}{2}}^{DS} = d_m + O(\Delta x),
\end{equation}
where the superscript $DS$ denotes the enhancement of the nonlinear weights \eqref{eq:omegas} using our novel method. 
We see, that neither the condition \eqref{eq:suff_cond}, nor \eqref{eq:cond3} is satisfied. 
However, as \eqref{eq:WENO_DS_JS_cond} holds, we can still guarantee convergence for the WENO-JS scheme with the smoothness indicators \eqref{eq:betann} with a convergence order degraded to the third order, cf.\ \cite{Borges}.

\subsection{Convergence analysis of WENO-Z scheme with new smoothness indicators $\beta_{m,i \pm \frac{1}{2}}^{DS}$}

Let us now analyse the convergence of the scheme \eqref{eq:WENOZ} with the new smoothness indicators \eqref{eq:betann}. 
From \eqref{eq:taylor_beta} we see that the smoothness indicators $\beta_{m,i \pm \frac{1}{2}}$ are of the form
\begin{equation}
\beta_{m,i \pm \frac{1}{2}} = f_{x}^2 \Delta x^2 + O(\Delta x^4) 
\end{equation}
and the global smoothness indicator \eqref{eq:tau}
\begin{equation}
\tau_5 = O(\Delta x^5).
\end{equation}
Then it holds
\begin{equation}
\begin{split}
   \beta_{m,i \pm \frac{1}{2}}^{DS} &= \beta_{m,i \pm \frac{1}{2}} (\delta_{m,i} + C) =  \bigl(f_{x}^2 \Delta x^2 + O(\Delta x^4) \bigr) \bigl(P(\bar{x}_{i}) + O(\Delta x)\bigr) \\
    &= f_{x}^2 P(\bar{x}_{i}) \Delta x^2 + O(\Delta x^3).
\end{split}\end{equation}
We take $\epsilon=0$ and choose $C$ such that $\Phi(\bar{x}_{i})+C>\kappa>0$, with $\kappa$ fixed. 
Then we see that in the non-critical points where $f_x \neq 0$
\begin{equation}
   \frac{\tau}{\beta_{m,i \pm \frac{1}{2}}^{DS}}= \hat{D}\Delta x^3 + O(\Delta x^4),
\end{equation}
where $\hat{D} = \frac{\frac{13}{3}|f_{xx}f_{xxx}|}{f_x^2 P(\bar{x}_{i})}$.
Substituting this into \eqref{eq:WENOZ} we obtain 
\begin{equation} 
  \alpha_{m,i \pm \frac{1}{2}}^{DS} = d_m\bigl(1+O(\Delta x^6)\bigr) \quad \text{ and } \quad
   \sum_{m=0}^{2}\alpha_{m,i \pm \frac{1}{2}}^{DS} = \bigl(1+O(\Delta x^6)\bigr),
\end{equation}
so it follows directly
\begin{equation}
   \omega_{m,i \pm \frac{1}{2}}^{DS} = d_m + O(\Delta x^6)
\end{equation}
and the condition \eqref{eq:suff_cond} is satisfied. 
Here we use the superscript $DS$ denoting the enhancement of the nonlinear weights \eqref{eq:WENOZ} using our novel method. 
 Since we ensure $P(\bar{x}_{i})>C>\kappa>0$, the multipliers $P(\bar{x}_{i})$ do not introduce any further critical points. 
Therefore the analysis of the critical points with $f_x = 0$ remains the same as in \cite{Borges}. 
Thus we can guarantee the fifth order convergence of the scheme \eqref{eq:WENOZ} with the smoothness indicators \eqref{eq:betann} also in the critical points.

\section{Numerical Results} \label{S:num}
For the system of ordinary differential equations resulting from \eqref{eq:hypsemi}
we use a third-order total variation diminishing (TVD) Runge-Kutta method \cite{Shu0} given by 
\begin{equation} \label{eq:runge_kutta}
\begin{split}
  u^{(1)} &= u^n + \Delta t L(u^n),  \\
  u^{(2)} &= \frac{3}{4}u^n + \frac{1}{4}u^{(1)} + \frac{1}{4}\Delta t L(u^{(1)}), \\
  u^{n+1} &= \frac{1}{3}u^n + \frac{2}{3}u^{(2)} + \frac{2}{3}\Delta t L(u^{(2)}), 
\end{split}\end{equation}
where $L = -\frac{1}{\Delta x}(\hat{f}_{i+\frac{1}{2}}-\hat{f}_{i-\frac{1}{2}})$ and $u^n$ is the solution at the time step $n$. 

For \eqref{eq:fluxsplit} we consider in our examples the Lax-Friedrichs flux splitting
\begin{equation} \label{eq:LF_flux_splitting}
    f^\pm(u) = \frac{1}{2} \bigl(f(u) \pm \alpha u\bigr),
\end{equation}
where $\alpha = \max\limits_u|f'(u)|$.

\subsection*{The Neural Network Structure}
The proposed neural network algorithm can be generally applied to any type of conservation laws. 
For the equations where discontinuities or shocks are present, we propose to train a neural network separately for each equation class. 
Then we can better adjust the size of a neural network and its structure as well as the loss function, which leads to better results.

As we mentioned earlier, the inputs to the CNN are the values $f^+(x_i)$ and $f^-(x_i)$, $i=0,\dots,N$ 
and we train two neural networks for a positive and negative part of a flux. 
(The superscripts $\pm$ will be further dropped.)

The first layer of the neural network is not learned, but represents a preprocessing of the solution from the last time step into a set of features that we assume to be suitable inputs for the following learned layers. 
Since our goal is to improve the smoothness indicators, we first compute the first and second central finite differences of $f(x_i)$, $i=0,\dots,N$. 
These parameters give us information about the smoothness of the solution and can facilitate and speed up the training of the CNN, and we can use a rather small CNN that still remains powerful. 
So we have the following values as input for the first learned hidden layer:
\begin{equation}\label{eq:preprocessing_1}
f_{\text{diff1}} = f(x_{i+1})-f(x_{i-1}), \quad f_{\text{diff2}} = f(x_{i+1})-2f(x_i)+f(x_{i-1}).
\end{equation}
The values $f_{\text{diff1}}$, $f_{\text{diff2}}$ computed from $f^+$ from \eqref{eq:LF_flux_splitting} represent the input values for the first neural network and the values $f_{\text{diff1}}$, $f_{\text{diff2}}$ computed from $f^-$ represent the input values for the second neural network.

Next, we use a fixed number of hidden layers, each with a specific kernel size and number of channels. 
We set these CNN parameters separately for each of the equation classes and experimentally find the best setting 
for each equation, keeping the size of proposed CNN small.
We move the kernel by one space step so the stride is set to 1, and we use an ELU activation function in all hidden layers except the last one where we use sigmoid.
In all our experiments, we set $C=0.1$ in \eqref{eq:betann}, which we experimentally found to be efficient. 
Let us note, that due to subsequent normalization of $\beta_m^{DS}$ during the computation of nonlinear weights, using large value of $C$ would decrease the effect of the trained multipliers. 
On the other hand, for $C$ close to zero the experimental order of convergence could be smaller on rough grids (but still achieved for $\Delta x \rightarrow 0$). We use the nonlinear weights as defined in \eqref{eq:WENOZ}, replacing $\beta_m$ with $\beta_m^{DS}$.
The value of $\epsilon$ is set to $10^{-13}$.

As the first choice of the loss function we use the mean square error 
\begin{equation} \label{eq:loss_L2}
   LOSS_{\rm MSE}(u) = \frac{1}{N} \sum_{i=0}^{N} (u_i - u_i^{\rm ref})^2, 
\end{equation}
where $u_i$ is a numerical approximation of $u(x_i)$ and $u_i^{\rm ref}$ is the corresponding reference solution. 
An advantage of this $L_2$-norm based loss function in contrast to the  $L_1$-norm based loss function are stronger gradients
with respect to $u_i$ resulting in faster training. 


\subsection{The Buckley-Leverett equation}
In the first example, we apply our neural network algorithm to the Buckley-Leverett equation, which was also considered, for example, in \cite{Shu0, Shu1, Jiang}. 
It is a typical example with a non-convex flux function modeling a two-phase fluid flow in a porous medium \cite{Leveque2002}. 
The flux is given by
\begin{equation} \label{eq:BL_flux}
    f(u) = \frac{u^2}{u^2 + a(1-u)^2}, \quad -1 \le x \le1, \quad 0\le t\le0.4,
\end{equation}
where $a<1$ is a constant indicating the ratio of the viscosities of the two fluids. 
The initial condition is set as
\begin{equation} \label{eq:BL_IC}
  u(x,0) = \begin{cases} 
    &1, \quad \text{if} \quad -0.5 \le x \le 0 , \\
   &0, \quad \text{if} \quad \text{elsewhere} \\
\end{cases} 
\end{equation}
and we use periodic boundary condition.

\begin{figure}[ht] 
\centering\includegraphics[width=1\linewidth]{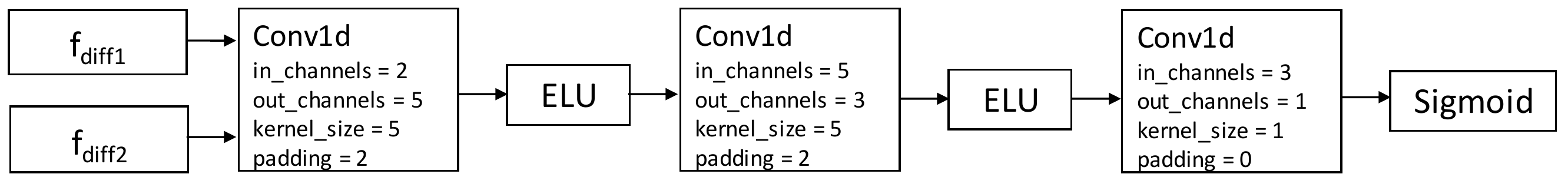}
\caption{A structure of the convolutional neural network used for the Buckley Leverett equation 
(structure is same for both inputs $f^+(x_i)$ and $f^-(x_i)$, $f_{\text{diff1}}$ and $f_{\text{diff2}}$ 
are defined in \eqref{eq:preprocessing_1} and are computed from both $f^+(x_i)$ and $f^-(x_i)$).} \label{fig:BL_conv}
\end{figure}

In our implementation, we use the CNN with 3 hidden layers with the structure described in the Figure~\ref{fig:BL_conv}. 
The inputs to the learned hidden layers are the features \eqref{eq:preprocessing_1}.  
First, we create the dataset for which we compute the reference solution for the equation \eqref{eq:hyp} with the flux \eqref{eq:BL_flux} and the initial condition \eqref{eq:BL_IC}. 
We randomly generate the parameter $a$ from a uniformly distributed range $[0.05,0.95]$. 
We divide the computational domain $[-1,1]$ into 1024 spatial steps and the solution is computed up to time $T=0.4$, 
where the time domain is divided into 8960 time steps. 
We use the WENO-Z method to compute this reference solution.

For training, we proceed as follows. At the beginning, we randomly choose a problem and its reference solution from our dataset. 
The weights of the CNN are randomly initialized and we train our model on a solution where our computational domain is divided into 128$\times$140 steps and successively compute the entire solution until the final time $T$. 
Using the solution at the time step $n$, we compute the solution at the time step $n+1$ and during this computation the CNN is used to predict the multipliers of the smoothness indicators.   
After each of these time steps, we compute the loss and its gradient with respect to the weights of the CNN using backpropagation algorithm. 
Then we use this gradient to update the weights, using the well-known Adam optimizer~\cite{Adam} with learning rate $0.0001$. 
After the last time step at time $T$, we test a model on a validation set and repeat the above steps. 
Then we select the model with the best performance on the validation set as our final model. 
For both training and comparing the performance of the models, we use the loss function defined as 
\begin{equation}
    LOSS(u) = LOSS_{\rm MSE}(u) + LOSS_{OF}(u),
\end{equation}
where $LOSS_{\rm MSE}(u)$ is defined in \eqref{eq:loss_L2} and 
\begin{equation}
    LOSS_{OF}(u) = \sum_{i=0}^{N} |\min(u_i, u_{\rm min}) - u_{\rm min}| + |\max(u_i, u_{\rm max}) - u_{\rm max}|
\end{equation}
represents the sum of the overflows of the solution above the maximum and below the minimum value of $u$, in our case $u_{\rm max}=1$ and $u_{\rm min}=0$. By adding this term to our loss function, we want to avoid the undesirable oscillations that occur especially in the first time steps of the solution. 

\begin{figure}[ht] 
\centering
\includegraphics[width=0.7\textwidth]{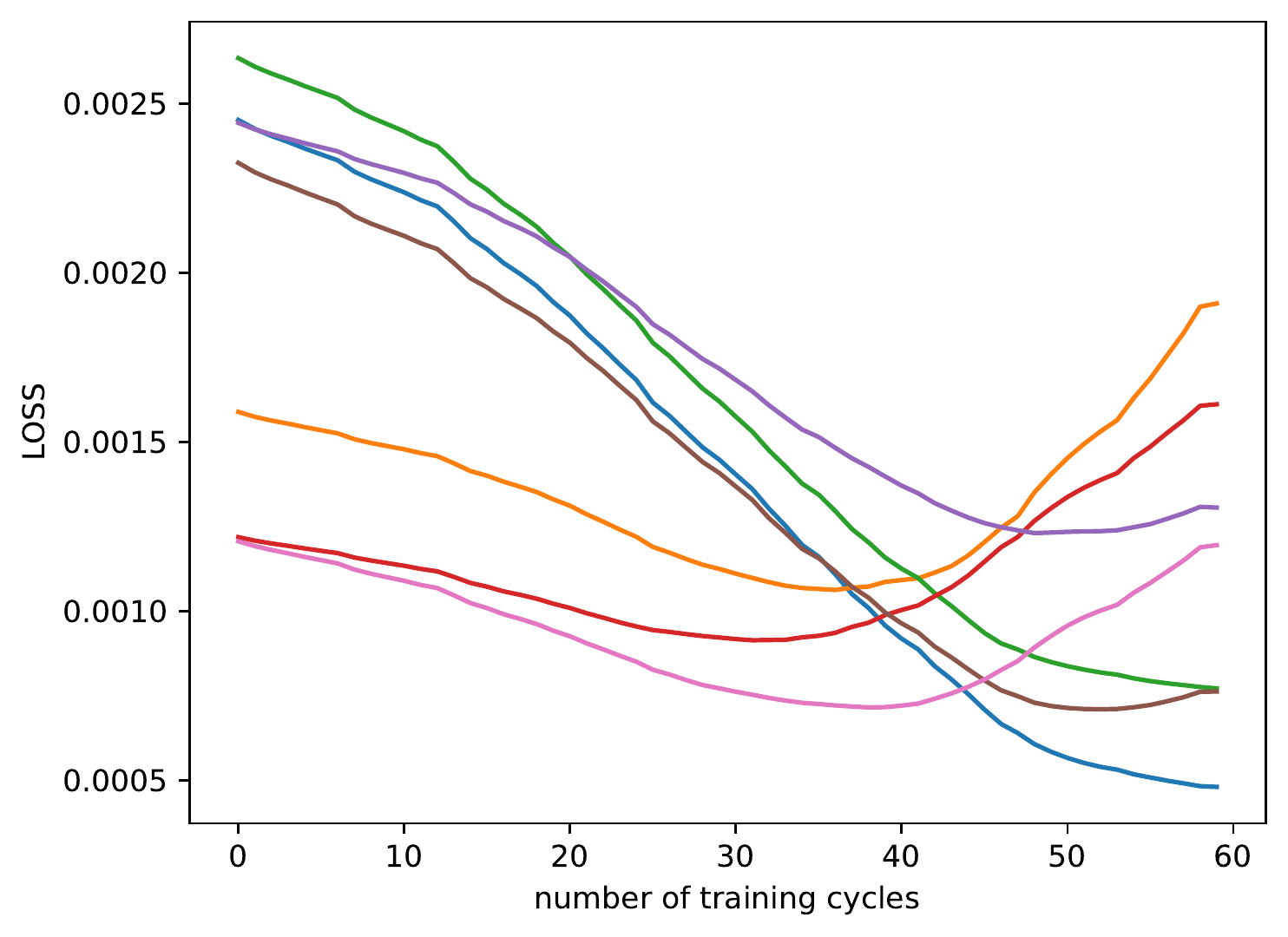}
\caption{Loss values for different validation problems at different training cycles ($x$-axis).}
\label{fig:BL_validation}
\end{figure}

The Figure~\ref{fig:BL_validation} shows how the value of the loss function for the problems from the validation set (which are not present in the training set) changes with increasing number of training cycles.
As training cycle we denote a sequence of training steps performed on a solution for a single randomly chosen parameter $a$ until the final time $T$. 
The loss is then evaluated at this final time $T$. 
We apparently see two optima for different values of $a$. 
If there are more than 50 training cycles, the loss begins to increase significantly for some problems, indicating that further training is not efficient.
This may be caused by overfitting or divergence, so an early stopping algorithm could be implemented efficiently. 
Since we want only one final numerical scheme as output, we choose the model obtained after the 46th training cycle as the final model and present the result computed with it.

We compare the $L_{\infty}$- and $L_2$-error in the Table~\ref{tab:BL} for the solution of the conservation law \eqref{eq:hyp} with \eqref{eq:BL_flux}, $a\in\{0.25, 0.4, 0.5, 0.6, 0.7, 0.8, 0.9 \}$. 
Let us note that these parameters were neither in the training, nor in the validation set.
\begin{table}[ht] 
\scalebox{0.75}{\begin{tabular}{|c|c|c|c|c|c|c|c|c|}
\hline
\multicolumn{1}{|c|}{}&\multicolumn{4}{|c|}{\ $L_\infty$}&\multicolumn{4}{|c|}{\ $L_2$} \\
\hline
$a$&  WENO-JS & WENO-Z &  WENO-DS & ratio &WENO-JS & WENO-Z &  WENO-DS & ratio\\
\hline
\hline
0.25 &  0.429654 &  0.435090 &  \textbf{0.183302} & 2.34 &  0.068405 &  0.067912 &  \textbf{0.034065} &1.99 \\  \hline
0.4 &  0.408252 &  0.405047 &  \textbf{0.340068} & 1.19& 0.059344 &  0.058160 & \textbf{ 0.056051} & 1.04 \\  \hline
0.5  &  0.317824 &  0.320094 &  \textbf{0.179696} &1.77 & 0.049913 &  0.049026 &  \textbf{0.033758} &1.45 \\  \hline
 0.6 &  0.459994 &  0.456687 &  \textbf{0.297523} &1.53  &0.062155 &  0.061275 &  \textbf{0.048766} &1 26 \\  \hline
0.7 &  0.476089 &  0.475015 &  \textbf{0.310196} &1.53  &0.073021 &  0.072581 &  \textbf{0.049836} &1.46 \\  \hline
0.8 &  0.207676 &  \textbf{0.197021} &  0.250032 & 0.79 &0.032560 &  \textbf{0.030994} &  0.038974 & 0.80\\  \hline
0.9 & 0.375720 &  0.367802 &  \textbf{0.181120} & 2.03 &0.062257 &  0.061834 &  \textbf{0.038510} & 1.61\\  \hline
\end{tabular}}
\caption{Comparison of $L_\infty$ and $L_2$ error of WENO-JS, WENO-Z and WENO-DS methods for the solution of the Buckley-Leverett equation with the initial condition \eqref{eq:BL_IC}. As 'ratio' we denote the minimum error of the methods WENO-JS and WENO-Z divided by the error of WENO-DS (rounded to 2 decimal points).}
\label{tab:BL}
\end{table}
We highlight the best performing WENO method using bold. In column denoted 'ratio' we divide the minimum error of WENO-JS and WENO-Z with an error of WENO-DS, showing how well our novel method works compared to the better one from the mentioned standard methods.
WENO-DS outperforms the standard WENO methods in most cases. 
For $a=0.8$, the error of WENO-DS is larger than the errors of the other two methods. 
However, this may be caused by the fact that the standard WENO methods perform disproportionately well for 
$a=0.8$ compared to other values of the parameter $a$.

\begin{figure}[ht] 
\begin{subfigure}{0.45\textwidth}
    \includegraphics[width=\textwidth]{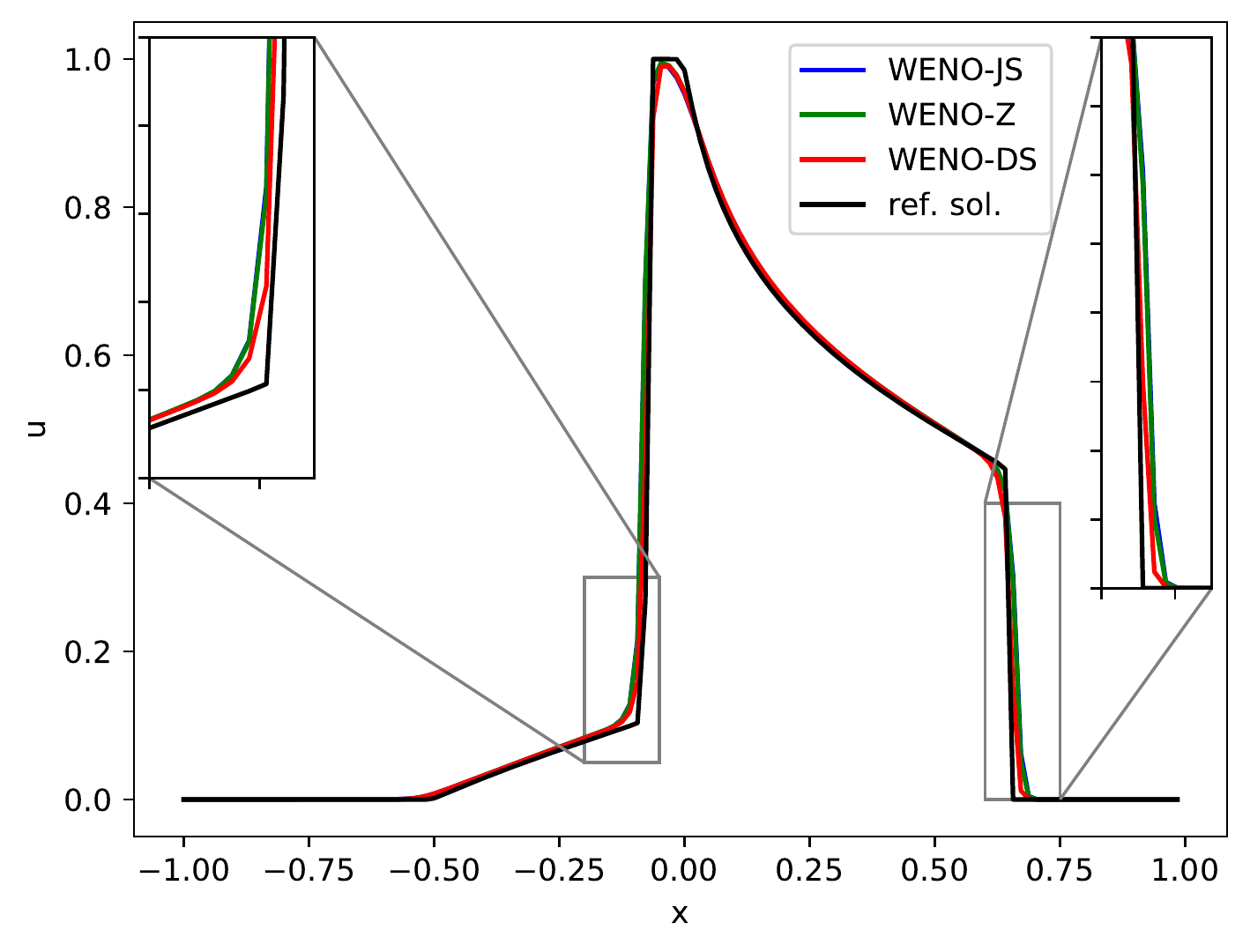}
    \caption{Solution for $a=0.25$}
\end{subfigure}
\begin{subfigure}{0.45\textwidth}
    \includegraphics[width=\textwidth]{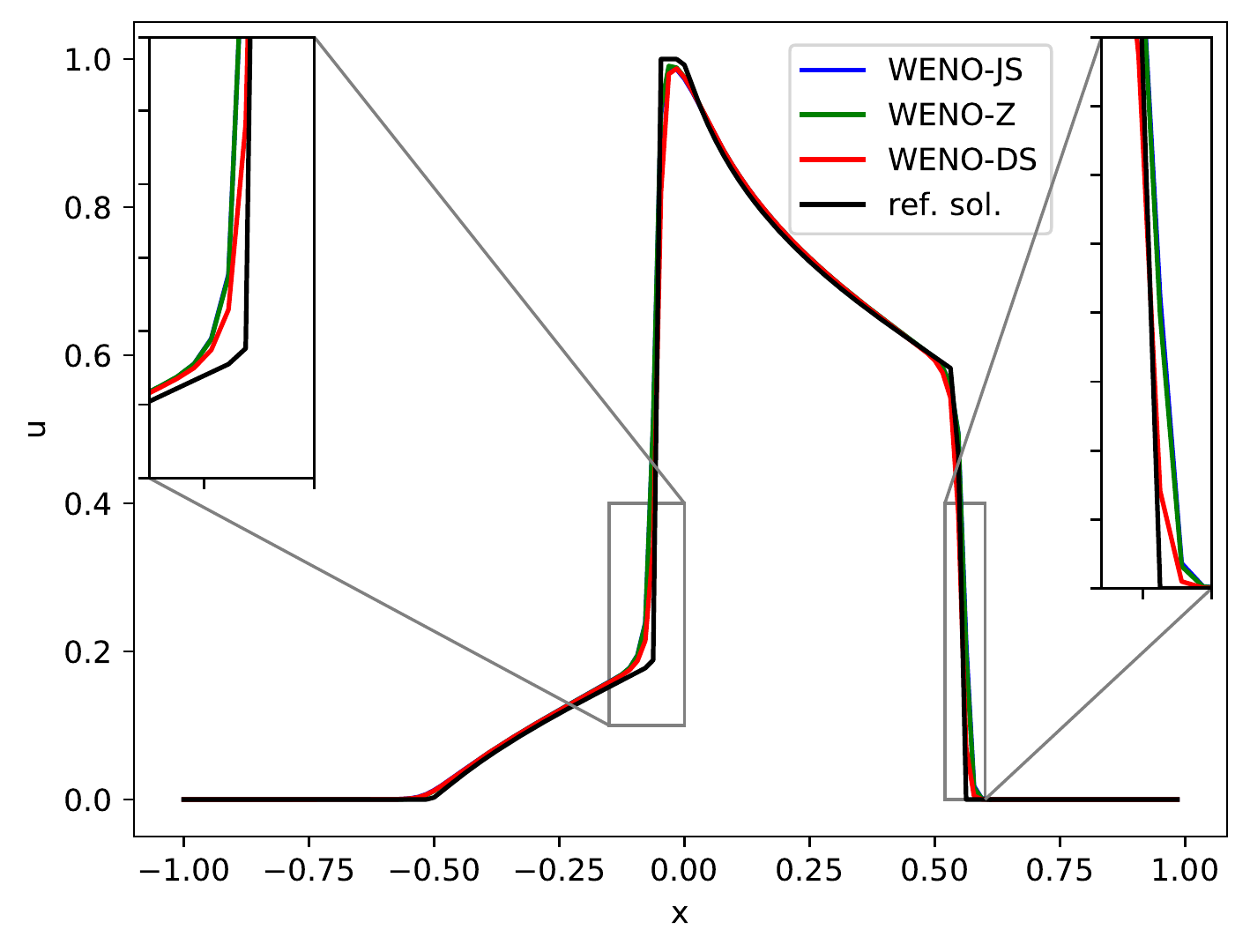}
    \caption{Solution for $a=0.5$}
\end{subfigure}
\caption{Comparison of the WENO-JS, WENO-Z and WENO-DS methods on the solution of 
the Buckley-Leverett equation with the initial condition \eqref{eq:BL_IC}, $N=128$.}
\label{fig:BL}
\end{figure}

In the Figure~\ref{fig:BL}, we present the solution of the Buckley-Leverett equation for 
the test problems with $a=0.25$ and $a=0.5$. 
We see that the WENO-DS gives a better solution quality than the WENO-JS or WENO-Z.

Finally, we verify the theoretically proven fifth-order convergence of the WENO-DS 
scheme for a transport equation with a smooth solution given as 
\begin{equation} \label{eq:transport}
     \frac{\partial u}{\partial t}+\frac{\partial u}{\partial x} = 0, \quad
u(x,0) =  sin(\pi x), \quad 0 \le x \le 2, \quad 0 \le t \le 0.5,
\end{equation}
with periodic boundary conditions. 
The results can be found in Table~\ref{tab:BL_order}. There is a great improvement when we compare our 
scheme with the WENO-NN scheme of Stevens and Colonius \cite{Stevens}, where the resulting scheme has only first-order convergence.
\begin{table}[ht]
\centering
\begin{tabular}{|c|c|c|c|c|}
\hline
\multicolumn{1}{|c|}{}&\multicolumn{2}{|c|}{\ WENO-Z}&\multicolumn{2}{|c|}{\ WENO-DS} \\
\hline
N &  $L_\infty$ & Order &  $L_\infty$ & Order    \\
\hline
\hline
20 &  9.369742e-03 &  - &  9.402549e-03 &  - \\ \hline
40 &  2.558719e-04 &  5.194516 &  2.558830e-04 &  5.199496 \\ \hline
80 &  9.466151e-06 &  4.756500 &  9.466165e-06 &  4.756560 \\ \hline
160 &  3.177833e-07 &  4.896663 &  3.177834e-07 &  4.896665 \\ \hline
320 &  9.957350e-09 &  4.996137 &  9.957351e-09 &  4.996138 \\ \hline
640 &  3.117835e-10 &  4.997145 &  3.117834e-10 &  4.997146 \\ \hline
\end{tabular}
\caption{$L_\infty$-norm error and convergence order of WENO-Z and WENO-DS on \eqref{eq:transport}.}
\label{tab:BL_order}
\end{table} 

\subsection{The inviscid Burgers’ equation}
In the next example we consider the inviscid Burgers' equation, where the flux function in \eqref{eq:hyp} is given by
\begin{equation}
    f(u) = \frac{u^2}{2}, \quad 0 \le x \le2, \quad 0 \le t \le 0.3. 
\end{equation}
We consider following initial conditions
\begin{gather} \label{eq:BE_IC_1}
  u(x,0) = \begin{cases} 
   &z_1, \quad \text{if} \quad 1 \le x \le 2 , \\
   &0, \quad \text{if} \quad \text{elsewhere}, \\
\end{cases} \\
u(x,0) = \exp (-z_2 (x - 1)^2), \label{eq:BE_IC_2} \\ 
u(x,0) = z_3 \sin (\pi x), \label{eq:BE_IC_3}
\end{gather}
where 
\begin{equation} \label{eq:z_range}
  z_1 \in \mathcal{U}[1,2], \quad z_2 \in \mathcal{U}[10,30], \quad z_3 \in \mathcal{U}[1,2].
\end{equation}
 Using these initial conditions, we cover problems with both continuous and discontinuous initial conditions, 
 and we simulate the shocks and discontinuities very well. We use periodic boundary condition.

\begin{figure}[ht] 
\centering\includegraphics[width=1\linewidth]{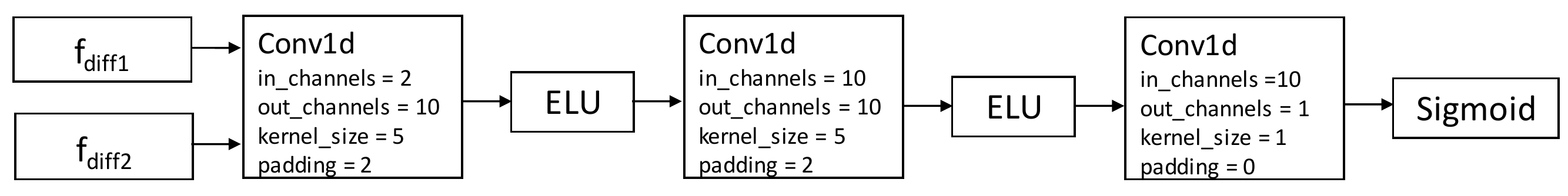}
\caption{The structure of the convolutional neural network used by the Burgers' equation 
(structure is same for both inputs $f^+(x_i)$ and $f^-(x_i)$, $f_{\text{diff1}}$ and $f_{\text{diff2}}$ 
are defined in \eqref{eq:preprocessing_1} and are computed from both $f^+(x_i)$ and $f^-(x_i)$).} \label{fig:BE_conv}
\end{figure}

We first create the data set for training, in which we compute the reference solution of the 
Burgers' equation with the initial conditions \eqref{eq:BE_IC_1}-\eqref{eq:BE_IC_3}. 
The computational domain is divided into 1024 space steps and 6400 time steps and the solution is computed up to time $T=0.3$ using the WENO-Z scheme. 
We use the CNN with 3 hidden layers with the structure described in the  Figure~\ref{fig:BE_conv}. 
For training, we proceed in the same way as in the previous example. 
The only differences are that the computational domain is divided into 128$\times$100 steps and the learning rate used by the Adam optimizer is $0.001$. 
In this example, we use the mean square error loss function \eqref{eq:loss_L2} for training and validation.
As the training on Burgers' equation exhibits much higher variance than in the Buckley-Leverett case,
we performed 3 trainings each with 90 training cycles and finally selected the model showing the best performance on the validation set.

\begin{table}[ht] 
\scalebox{0.68}{ \begin{tabular}{|p{1.7cm}|c|c|c|c|c|c|c|c|c|}
\hline
\multicolumn{1}{|c|}{}&\multicolumn{1}{|c|}{}&\multicolumn{4}{|c|}{\ $L_\infty$}&\multicolumn{4}{|c|}{\ $L_2$} \\
\hline
\centering initial  condition &  \multirow{2}*{$z_j$} &  \multirow{2}*{WENO-JS} & \multirow{2}*{WENO-Z} & \multirow{2}*{WENO-DS} & \multirow{2}*{ratio}&\multirow{2}*{WENO-JS} & \multirow{2}*{WENO-Z} &  \multirow{2}*{WENO-DS} & \multirow{2}*{ratio}\\
\hline
\hline
\centering  \eqref{eq:BE_IC_1} & 1.19  &  0.632213 &  0.632788 &  \textbf{0.303060} & 2.09&  0.082373 &  0.082141 &  \textbf{0.046932} &1.75\\ \hline
& 1.53 &  0.594943 &  0.58678 &  \textbf{0.485714} & 1.21&  0.080341 &  0.07877 &  \textbf{0.067175} & 1.17\\ \hline
 & 1.84 &  0.704680 &  0.694358 &  \textbf{0.542599} &1.28 &  0.094967 &  0.093019 &  \textbf{0.076102} & 1.22\\ \hline
\centering \eqref{eq:BE_IC_2} & 14.94 &  0.113498 &  0.104374 &  \textbf{0.100061} &1.04 &  0.016926 & \textbf{0.015137} &  0.015164  &1.00\\ \hline
 & 21.65 &  0.236125 &  0.229290 &  \textbf{0.196110} &1.17  & 0.032979 &  0.031680 &  \textbf{0.029141} & 1.09\\ \hline
 & 29.08 &  0.312595 &  \textbf{0.310937} &  0.388739 & 0.80 & 0.040632 &  \textbf{0.040199} &  0.049385 & 0.81\\ \hline
 \centering \eqref{eq:BE_IC_3} & 1.46 &  0.059072 &  0.056751 &  \textbf{0.051553} & 1.10&  0.010443 &  0.010032 &  \textbf{0.007307} & 1.37\\ \hline
 & 1.6 &  0.063780 &  0.061391 &  \textbf{0.037165} &1.65 &  0.011275 &  0.010853 &  \textbf{0.005552} & 1.95\\ \hline
 & 1.9 &  0.072586 &  0.069995 &  \textbf{0.023841} &2.94 & 0.012831 &  0.012373 &  \textbf{0.003396} & 3.64\\ \hline
\end{tabular}}
\caption{Comparison of $L_\infty$ and $L_2$ error of WENO-JS, WENO-Z and WENO-DS methods for the solution of the Burgers' equation with the initial condition parameters inside of training set intervals \eqref{eq:z_range}. As 'ratio' we denote the minimum error of the methods WENO-JS and WENO-Z divided by the error of WENO-DS (rounded to 2 decimal points).}
\label{tab:BE_1}
\end{table}

\begin{table}[ht] 
\scalebox{0.68}{ \begin{tabular}{|p{1.7cm}|c|c|c|c|c|c|c|c|c|}
\hline
\multicolumn{1}{|c|}{}&\multicolumn{1}{|c|}{}&\multicolumn{4}{|c|}{\ $L_\infty$}&\multicolumn{4}{|c|}{\ $L_2$} \\
\hline
\centering initial  condition &  \multirow{2}*{$z_j$} &  \multirow{2}*{WENO-JS} & \multirow{2}*{WENO-Z} & \multirow{2}*{WENO-DS} & \multirow{2}*{ratio}&\multirow{2}*{WENO-JS} & \multirow{2}*{WENO-Z} &  \multirow{2}*{WENO-DS} & \multirow{2}*{ratio}\\
\hline
\hline
 \centering  \eqref{eq:BE_IC_1} & 0.71 &  0.204162 &  0.199246 &  \textbf{0.150400} &1.32 &  0.031740 &  0.030740 &  \textbf{0.028178} &1.09 \\ \hline
 & 2.41 &  1.541714 &  1.554201 &  \textbf{1.004904} & 1.53 & 0.199285 &  0.200194 &  \textbf{0.129223} &1.54 \\ \hline
 & 2.57 &  1.063823 &  1.055948 &  \textbf{0.755908} & 1.40 & 0.140068 &  0.138411 &  \textbf{0.102963} &1.34 \\ \hline
 & 3.13 &  0.622858 &  0.600619 &  \textbf{0.287540} & 2.09 & 0.087067 &  0.084495 &  \textbf{0.054477} &1.55 \\ \hline
 \centering \eqref{eq:BE_IC_2} & 33.9 &  0.351086 &  0.345278 &  \textbf{0.266422} & 1.30&  0.045665 &  0.044691 &  \textbf{0.036426} & 1.23\\ \hline
 & 34.67 &  0.285791 &  0.283237 &  \textbf{0.194424} &1.46  & 0.037350 &  0.036815 &  \textbf{0.027150} &1.36 \\ \hline
 \centering \eqref{eq:BE_IC_3} & 0.94 &  0.009524 &  \textbf{0.007189} &  0.007898 & 0.91 & 0.001737 &  0.001509 &  \textbf{0.001491} & 1.01\\ \hline
 & 2.12 &  0.077503 &  0.074744 &  \textbf{0.012296} & 6.08&  0.013701 &  0.013213 &  \textbf{0.001712} & 7.72\\ \hline
 & 2.44 &   0.083010 &  0.080022 &  \textbf{0.003978} & 20.12&  0.014675 &  0.014147 &  \textbf{0.000537} & 26.35\\ \hline
\end{tabular}}
\caption{Comparison of $L_\infty$ and $L_2$ error of WENO-JS, WENO-Z and WENO-DS methods for the solution of the Burgers' equation with the initial condition parameters outside of training set intervals \eqref{eq:z_range}. As 'ratio' we denote the minimum error of the methods WENO-JS and WENO-Z divided by the error of WENO-DS (rounded to 2 decimal points).}
\label{tab:BE_2}
\end{table}

\begin{figure}[ht!]
    \centering
    \begin{subfigure}[b]{0.45\textwidth}
            \includegraphics[width=\textwidth]{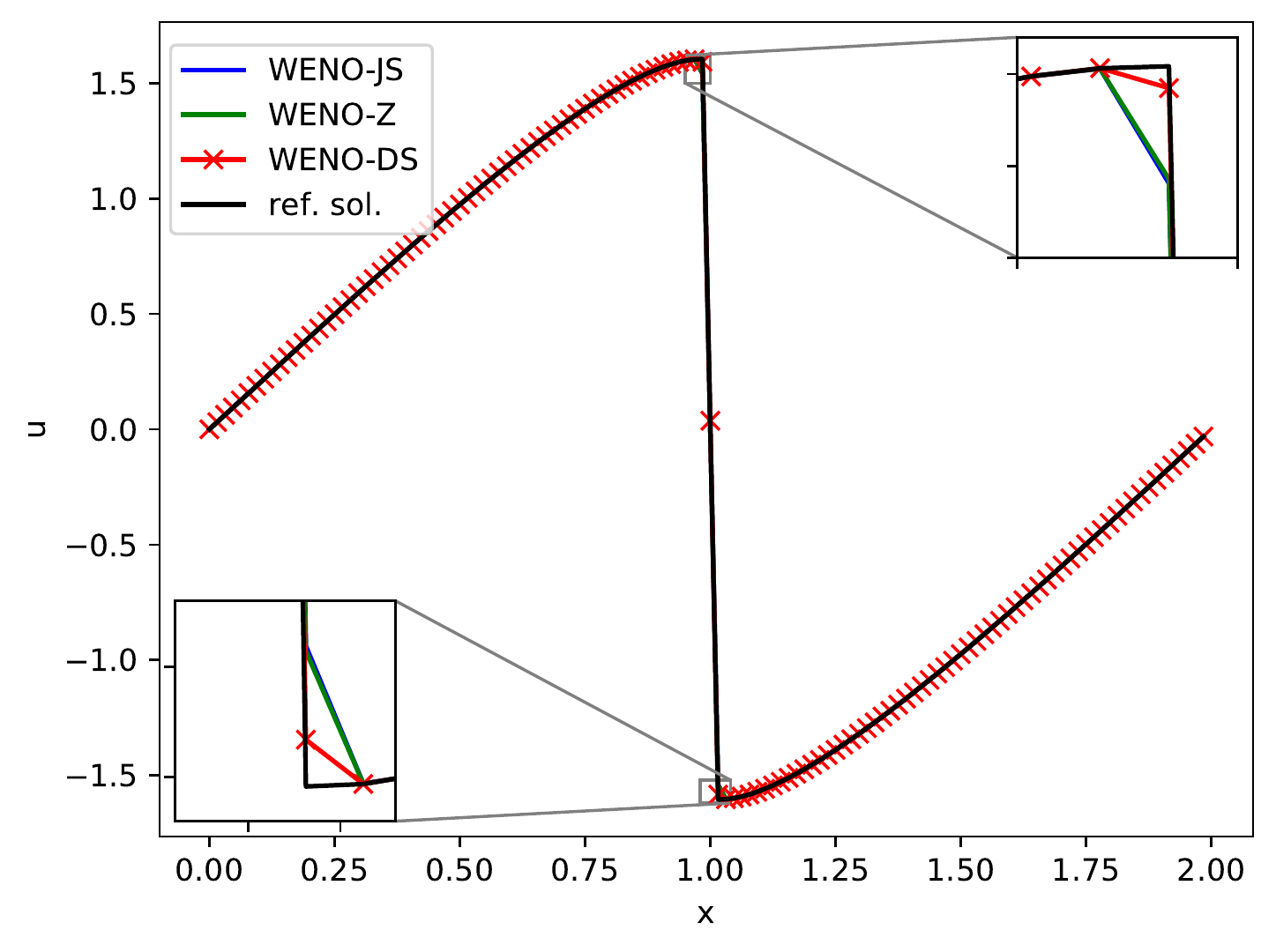}
            \caption{Initial condition \eqref{eq:BE_IC_3} with $z_3 = 1.6$.}
            \label{fig:BE_a}
    \end{subfigure}
    \begin{subfigure}[b]{0.45\textwidth}
            \includegraphics[width=\textwidth]{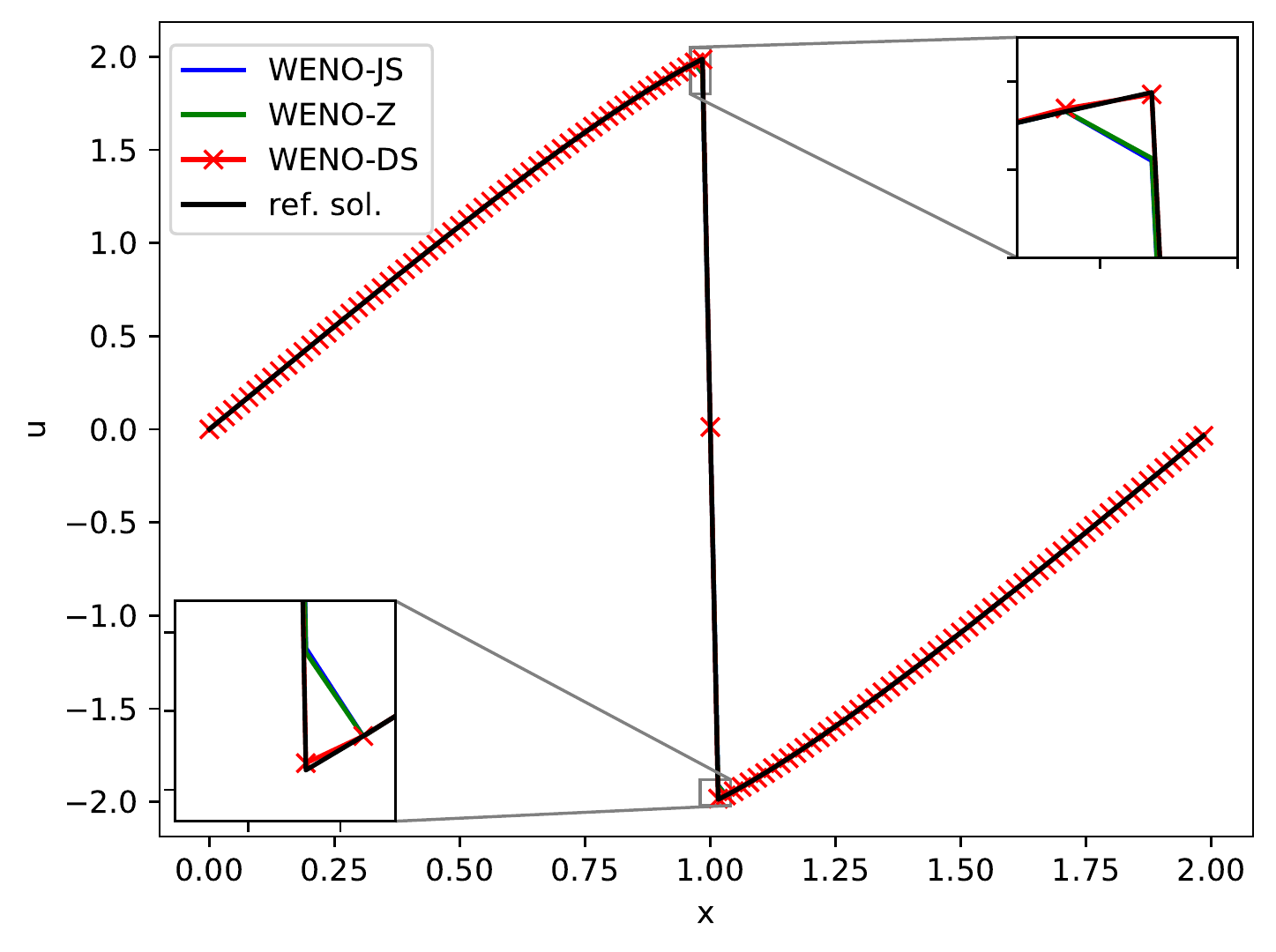}
            \caption{Initial condition \eqref{eq:BE_IC_3} with $z_3 = 2.12$.}
            \label{fig:BE_b}
    \end{subfigure}
    \begin{subfigure}[b]{0.45\textwidth}
            \includegraphics[width=\textwidth]{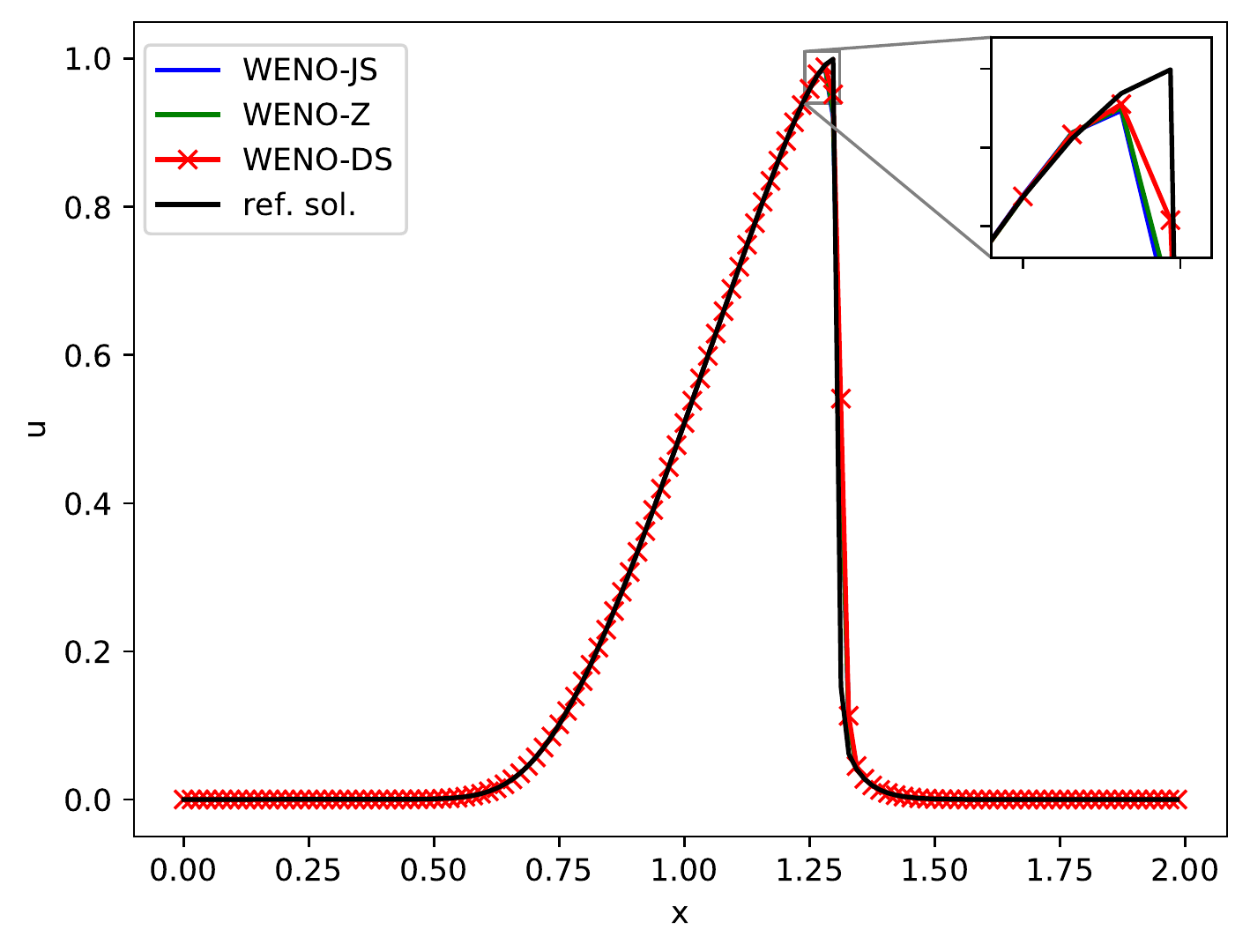}
            \caption{Initial condition \eqref{eq:BE_IC_2} with $z_2 = 29.08$.}
            \label{fig:BE_c}
    \end{subfigure}
    \begin{subfigure}[b]{0.45\textwidth}
            \includegraphics[width=\textwidth]{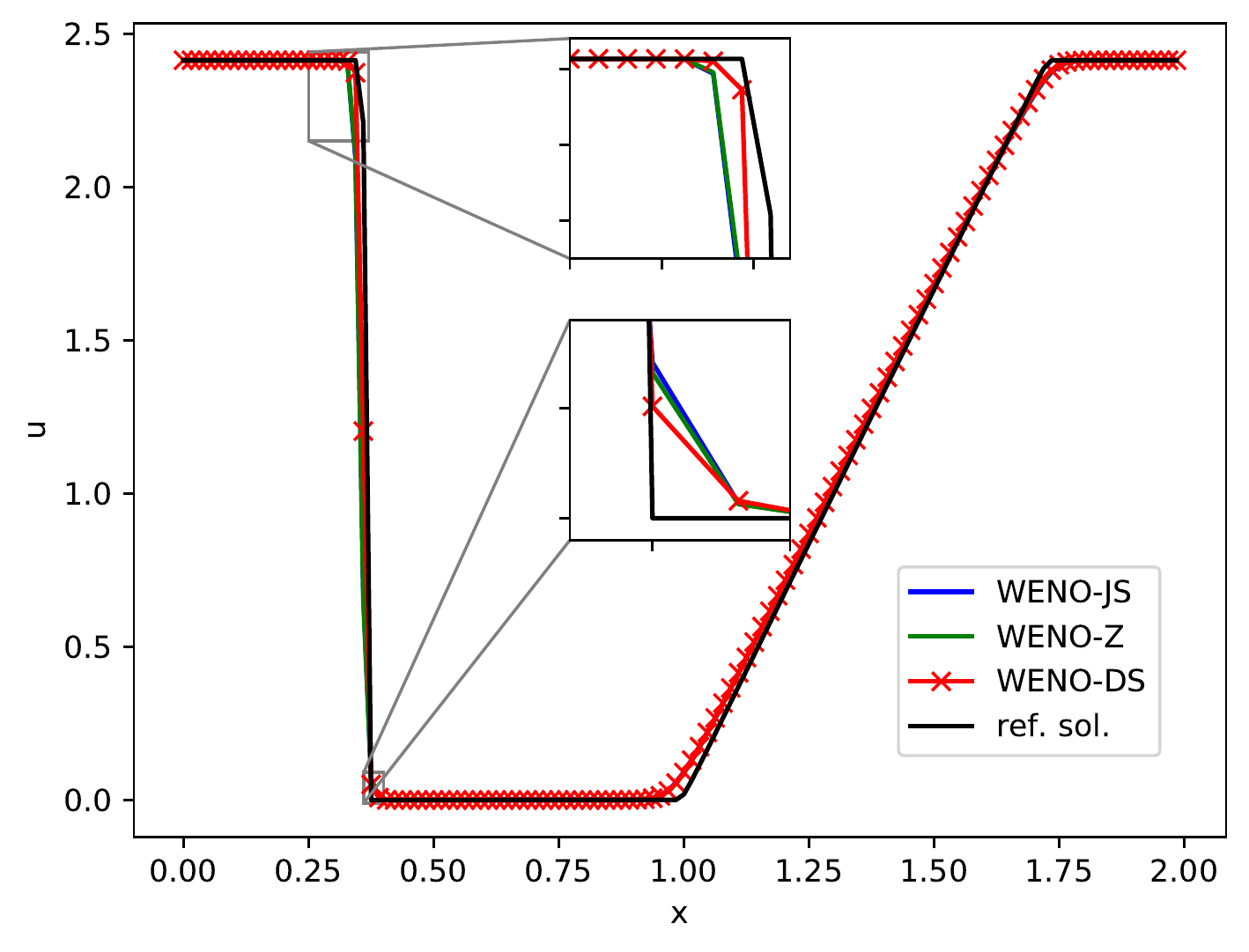}
            \caption{Initial condition \eqref{eq:BE_IC_1} with $z_1 = 2.41$.}
            \label{fig:BE_d}
    \end{subfigure}
    \caption{Comparison of the WENO-JS, WENO-Z and WENO-DS methods for the solution of the Burgers' equation with various initial conditions, $N=128$.}\label{fig:BE}
\end{figure}

We compare the errors on the problems from the test set in Table~\ref{tab:BE_1} and \ref{tab:BE_2}. 
These were not in the training or validation set and the parameters were randomly generated.
We observe rather small or no improvement for problems with the initial condition \eqref{eq:BE_IC_2}, but the improvement is significant for the solution with the discontinuous initial condition \eqref{eq:BE_IC_1} as well as with the initial condition \eqref{eq:BE_IC_3}.

We conclude that the WENO-DS significantly outperforms the classical WENO methods. 
It should be noted that although our training set was created with the parameters sampled from uniform distribution as specified in \eqref{eq:z_range}, the method can also generalise for parameter values outside of these intervals, as can be seen in Table~\ref{tab:BE_2}. 
Especially, we highlight the last two problems from Table~\ref{tab:BE_2}, where we see a great improvement.

In the Figure~\ref{fig:BE} we show the solution of the Burgers' equation with the initial condition \eqref{eq:BE_IC_1} for $z_1=2.41$, \eqref{eq:BE_IC_2} for $z_2=29.08$, \eqref{eq:BE_IC_3} for $z_3=1.6$ and \eqref{eq:BE_IC_3} for $z_3=2.12$. 
We observe that WENO-DS captures shocks and discontinuities very well and gives us a better solution compared to WENO-JS and WENO-Z.

\subsection{The one-dimensional Euler equations} 
We now investigate how WENO-DS behaves when applied to the one-dimensional Euler system, which is considered a classical benchmark problem for methods for conservation laws. 
It has the form
\begin{equation}
    \begin{split}
        \frac{\partial \rho}{\partial t} + \frac{\partial(\rho u)}{\partial x} &= 0, \\
    \frac{\partial \rho u}{\partial t} + \frac{\partial(\rho u^2 + p)}{\partial x} &= 0, \\
    \frac{\partial E}{\partial t} + \frac{\partial(u E + u p)}{\partial x} &= 0, 
    \end{split}
\end{equation}
where $\rho$ is the density, $u$ is the velocity, 
$p$ is the pressure and $E$ is a total energy given by
\begin{equation}
    E = \frac{p}{\gamma-1} + \frac{1}{2}\rho u^2.
\end{equation}
We take $\gamma=1.4$, which is the ratio of the specific heats. 
To compute the fluxes, we use the characteristic decomposition of the system according to the steps in \cite{Shu1998}. 
We use the Roe scheme to obtain the eigenvectors and eigenvalues \cite{Roe} and the Lax-Friedrichs flux splitting to obtain the corresponding component of the flux. 
We take the solution based on \cite{Wesseling} as the reference solution. 

One of the most common benchmark problems is the Sod problem \cite{Sod}, where the initial condition is specified as
\begin{equation}  \label{eq:sod}
  (\rho,u,p) = \begin{cases} 
    (1, 0, 1) \quad 0 \le x \le 0.5,\\
    (0.125, 0, 0.1) \quad 0.5 < x \le 1 \\
\end{cases} 
\end{equation}
and the solution is computed up to time $T = 0.1$. 
We use an adaptive step size
\begin{equation}
    \Delta t = \frac{0.9 \Delta x}{\max(c_i+|u_i|)}, \quad c^2= \frac{\gamma p}{\rho},
\end{equation}
where $u_i$ is the local velocity and $c_i$ the local speed of sound.

The solution consists of the left rarefaction wave, the right travelling contact wave and the right shock wave. 
We want to imitate this behavior of the solution, so we construct our data set as described in \ref{A:1}.

\begin{figure}[ht]
\centering\includegraphics[width=1\linewidth]{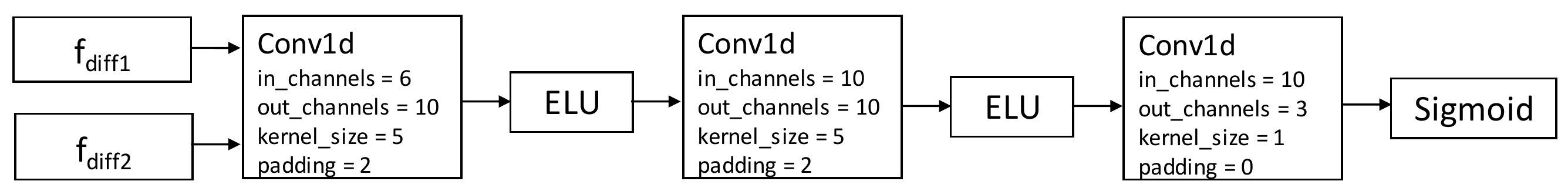}
\caption{A structure of the convolutional neural network used in the Euler system 
(structure is same for both inputs $f^+(x_i)$ and $f^-(x_i)$, $f_{\text{diff1}}$ and $f_{\text{diff2}}$ are defined in \eqref{eq:preprocessing_1} and are computed from both $f^+(x_i)$ and $f^-(x_i)$).} 
\label{fig:Euler_conv}
\end{figure}

We use the CNN with 3 hidden layers with the structure described in Figure~\ref{fig:Euler_conv}. 
After projecting the flux and the solution on the characteristic fields using the left eigenvectors, 
we use Lax-Friedrichs flux splitting for each component of characteristic variables. 
From these values we compute the features \eqref{eq:preprocessing_1}, which are the inputs to the learned hidden layers.

\begin{figure}[b!]
\begin{subfigure}{.32\textwidth}
\centering
\includegraphics[width=\linewidth]{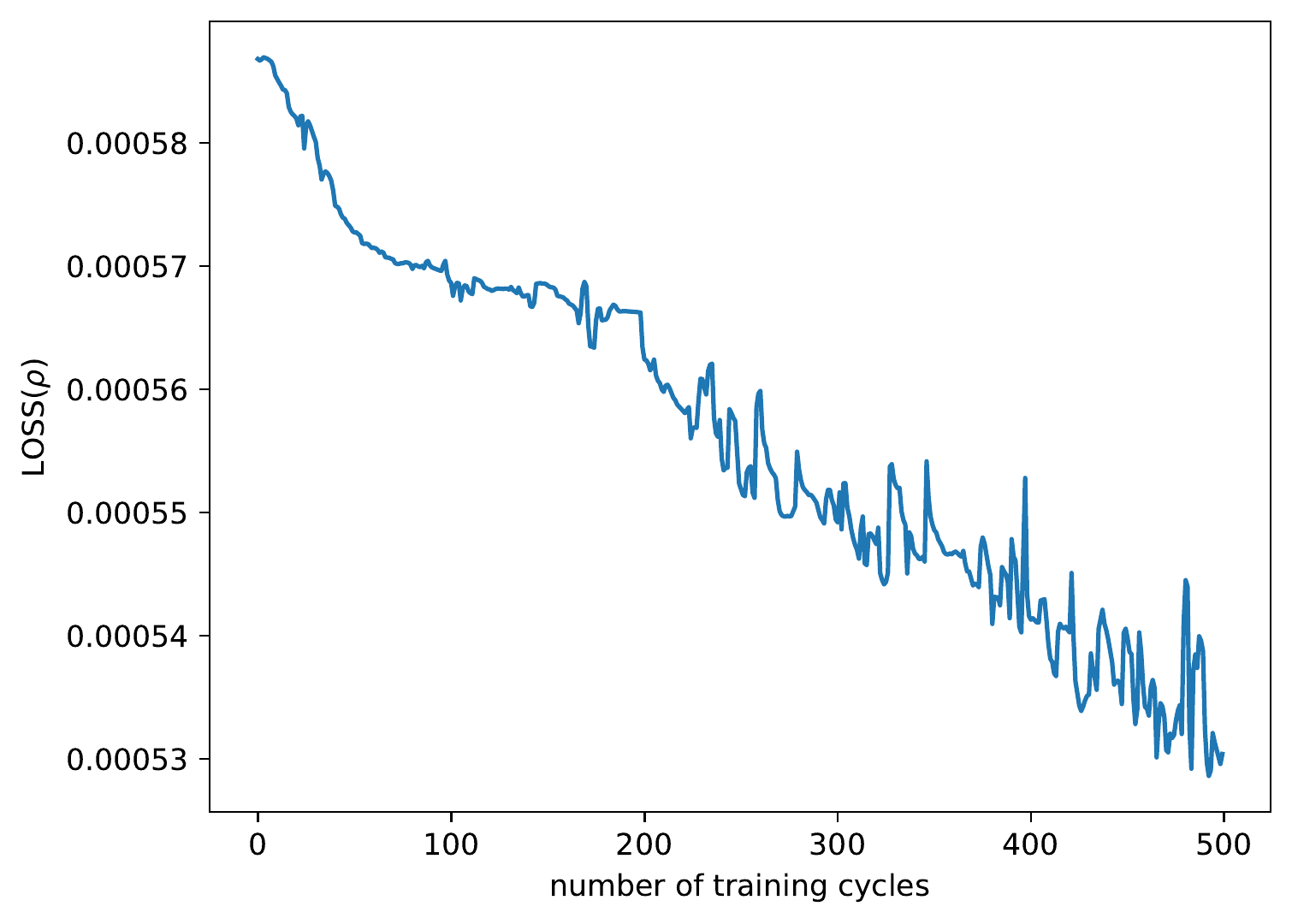} \label{fig:loss_rho}
\caption{$LOSS_{\rm MSE}(\rho)$}
\end{subfigure}
\begin{subfigure}{.32\textwidth}
\centering
\includegraphics[width=\linewidth]{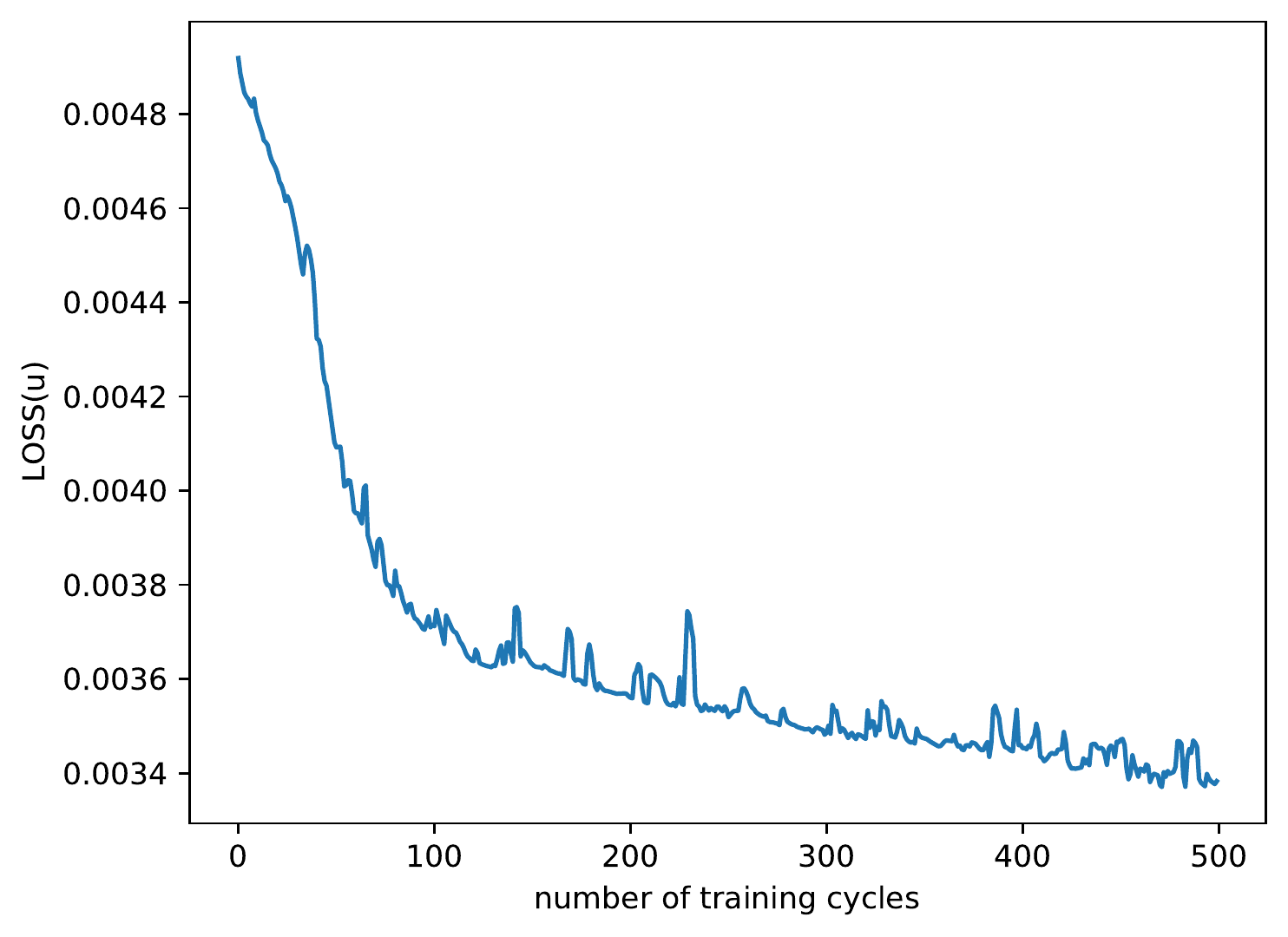} \label{fig:loss_u}
\caption{$LOSS_{\rm MSE}(u)$}
\end{subfigure}
\begin{subfigure}{.32\textwidth}
\centering
\includegraphics[width=\linewidth]{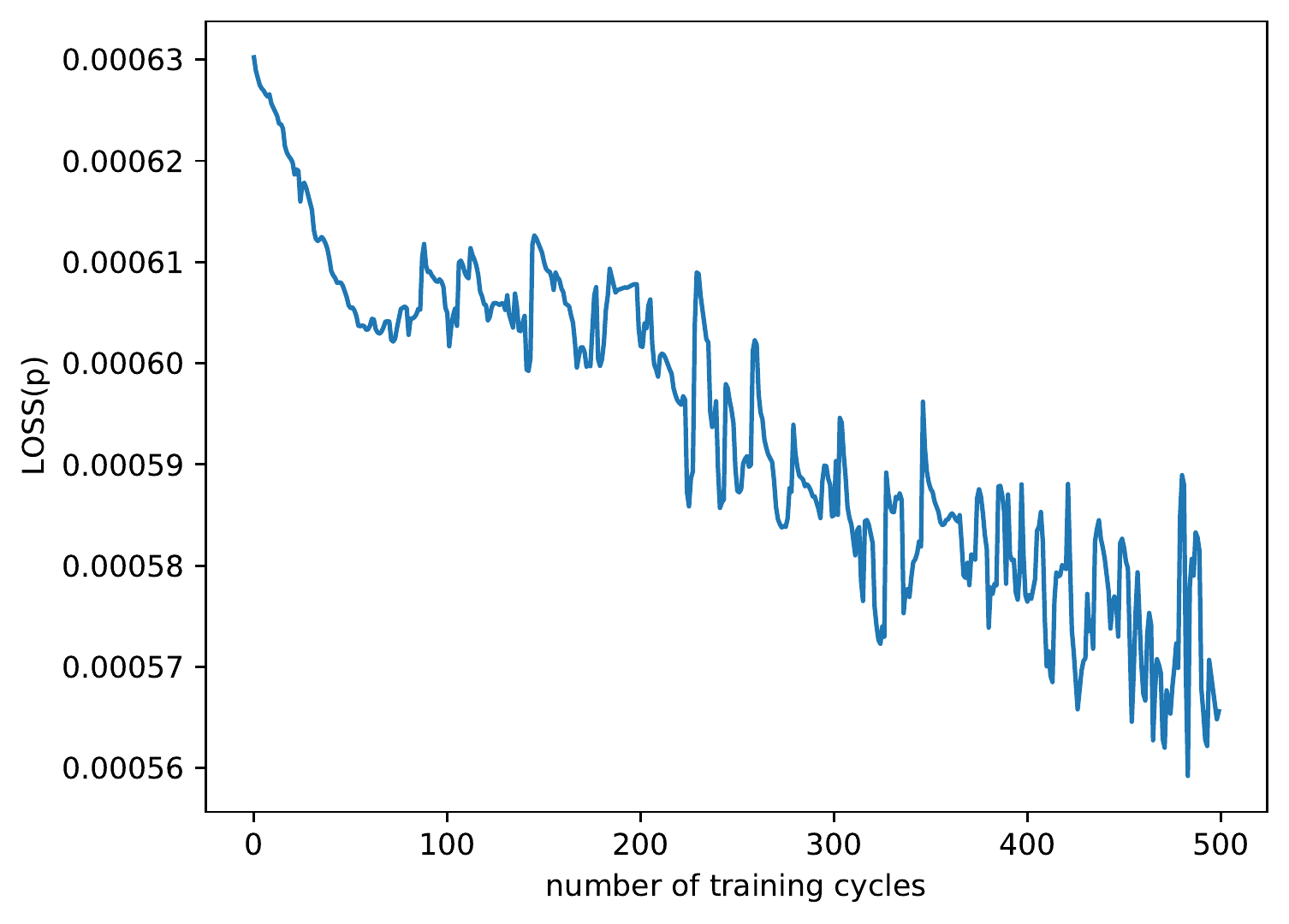} \label{fig:loss_p}
\caption{$LOSS_{\rm MSE}(p)$}
\end{subfigure}
\caption{Development of $LOSS_{\rm MSE}(\cdot)$ values with increasing number of training cycles.} \label{fig:Euler_loss}
\end{figure}
In this example we repeat the training procedure of the previous examples with some small modifications. 
To begin, we randomly generate the initial state from the dataset described earlier. 
We divide the spatial domain into 64 steps and compute the solution for the given initial state up to the time $T=0.1$. 
After each time step we compute loss using the reference solution from \cite{Wesseling}. 
We use the gradient to update the weights, using Adam optimizer with learning rate $0.001$. 
At the last time step we test the model and repeat the procedure with the new initial parameters $(\rho, u, p)$.
We use the loss function
\begin{equation} \label{eq:Euler_loss}
    LOSS(\rho, u, p) = LOSS_{\rm MSE}(\rho) + LOSS_{\rm MSE}(u) + LOSS_{\rm MSE}(p)
\end{equation}
for training and validation.

We show in the Figure \ref{fig:Euler_loss} how the values of $LOSS_{\rm MSE}(\cdot)$ from \eqref{eq:Euler_loss} develop with the increasing number of training cycles. 
These values were obtained by testing the method on the validation problem, which was the Sod problem with the initial condition \eqref{eq:sod}.

\begin{figure}[htp!]
\begin{subfigure}{.32\textwidth}
\centering
\includegraphics[width=\linewidth]{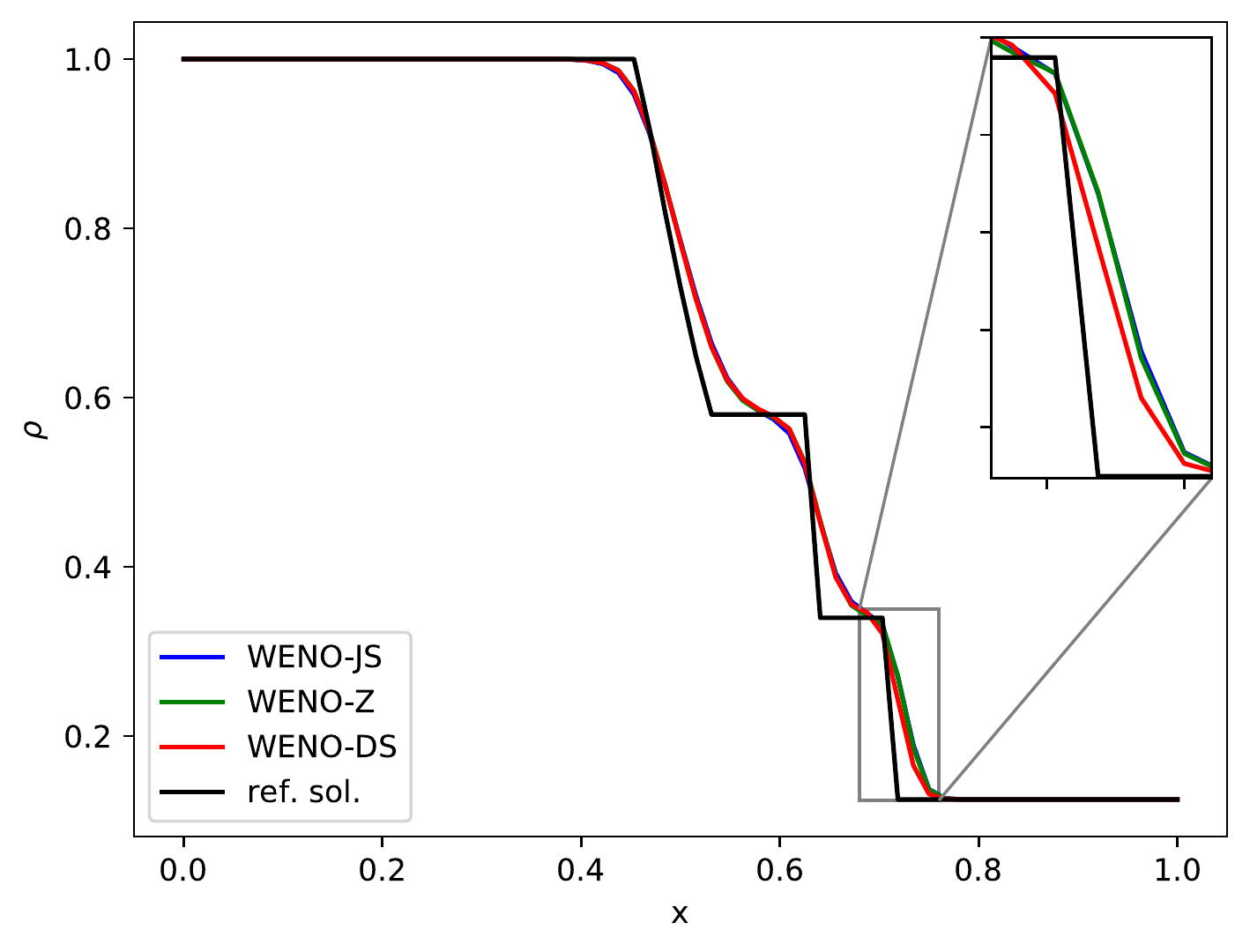} \label{fig:Sod_rho_72}
\end{subfigure}
\begin{subfigure}{.32\textwidth}
\centering
\includegraphics[width=\linewidth]{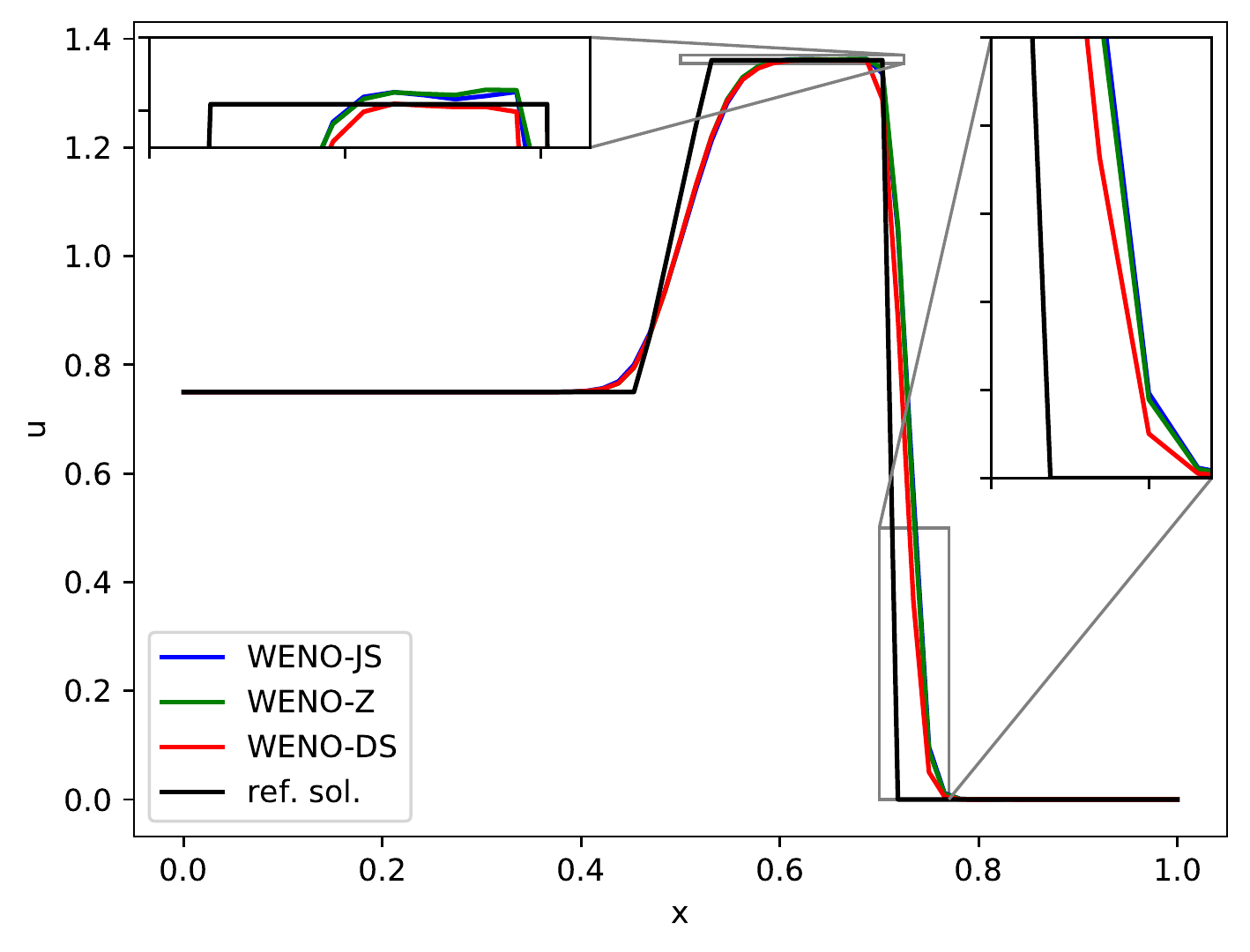} \label{fig:Sod_u_72}
\end{subfigure}
\begin{subfigure}{.32\textwidth}
\centering
\includegraphics[width=\linewidth]{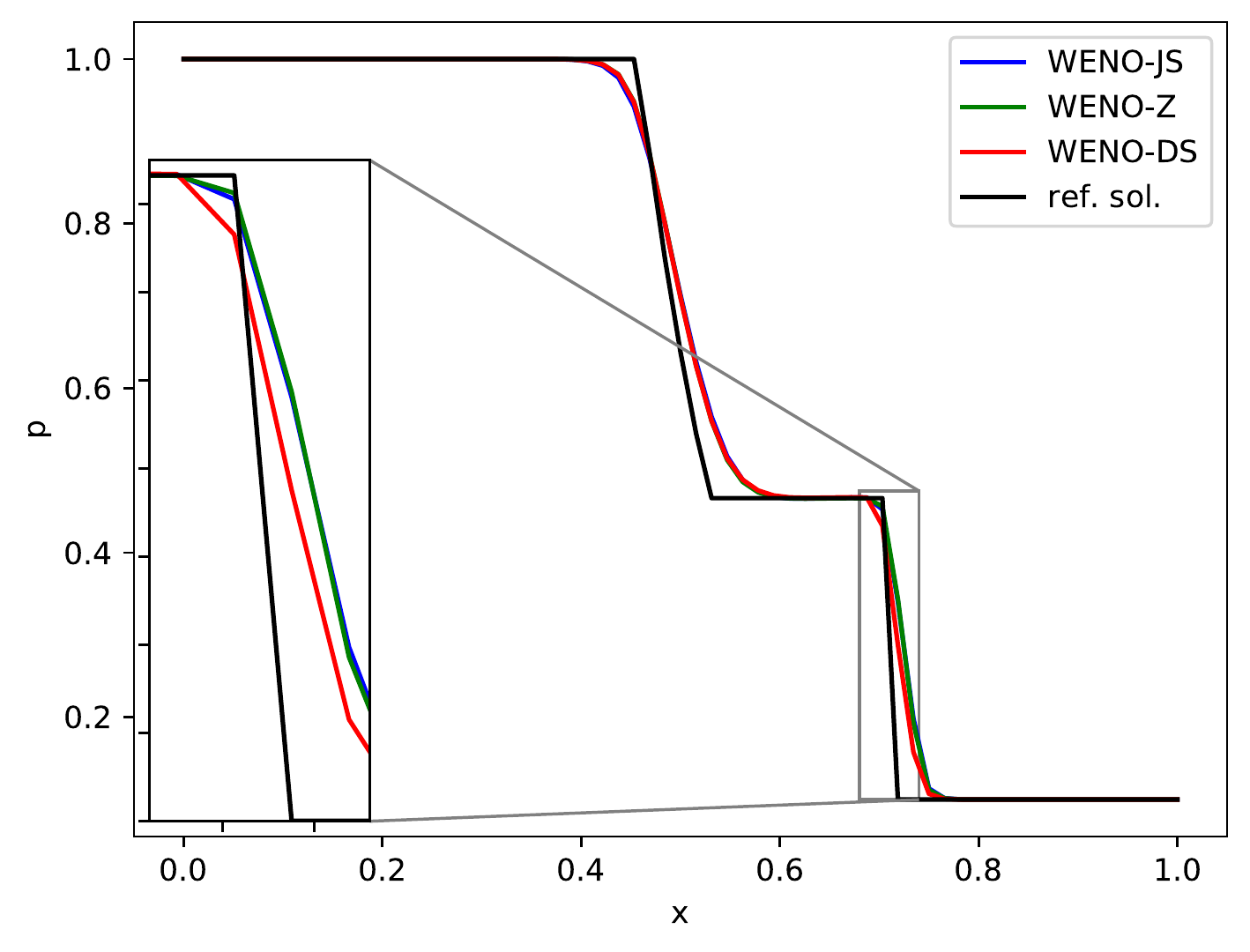} \label{fig:Sod_p_72}
\end{subfigure}
\caption{Solution of Sod problem \eqref{eq:sod_mod}, using the WENO-JS, WENO-Z and WENO-DS obtained after the 72th training cycle, $N = 64$.} \label{fig:Sod_72}
\end{figure}

\begin{figure}[htp!]
\begin{subfigure}{.32\textwidth}
\centering
\includegraphics[width=\linewidth]{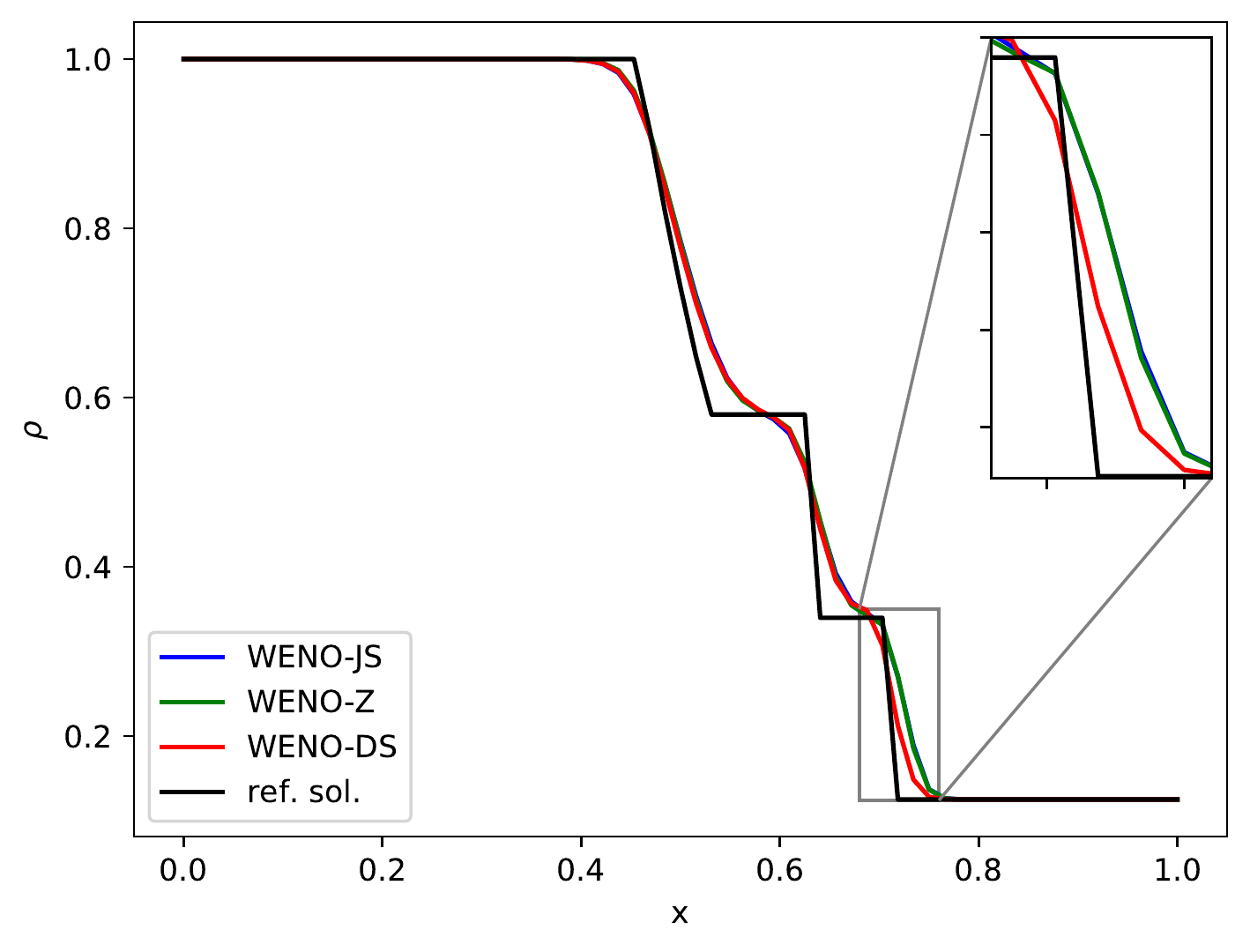} \label{fig:Sod_rho_483}
\end{subfigure}
\begin{subfigure}{.32\textwidth}
\centering
\includegraphics[width=\linewidth]{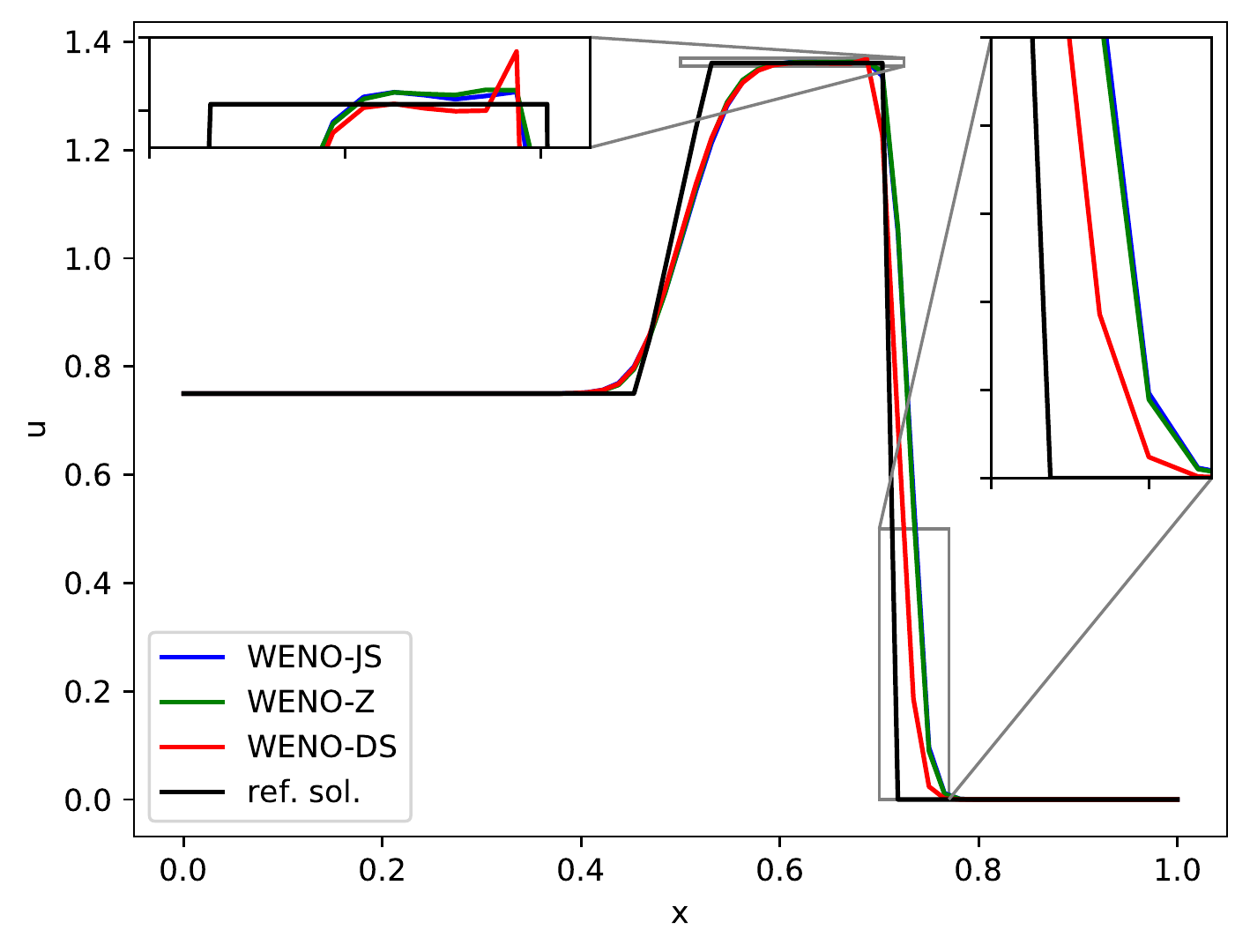} \label{fig:Sod_u_483}
\end{subfigure}
\begin{subfigure}{.32\textwidth}
\centering
\includegraphics[width=\linewidth]{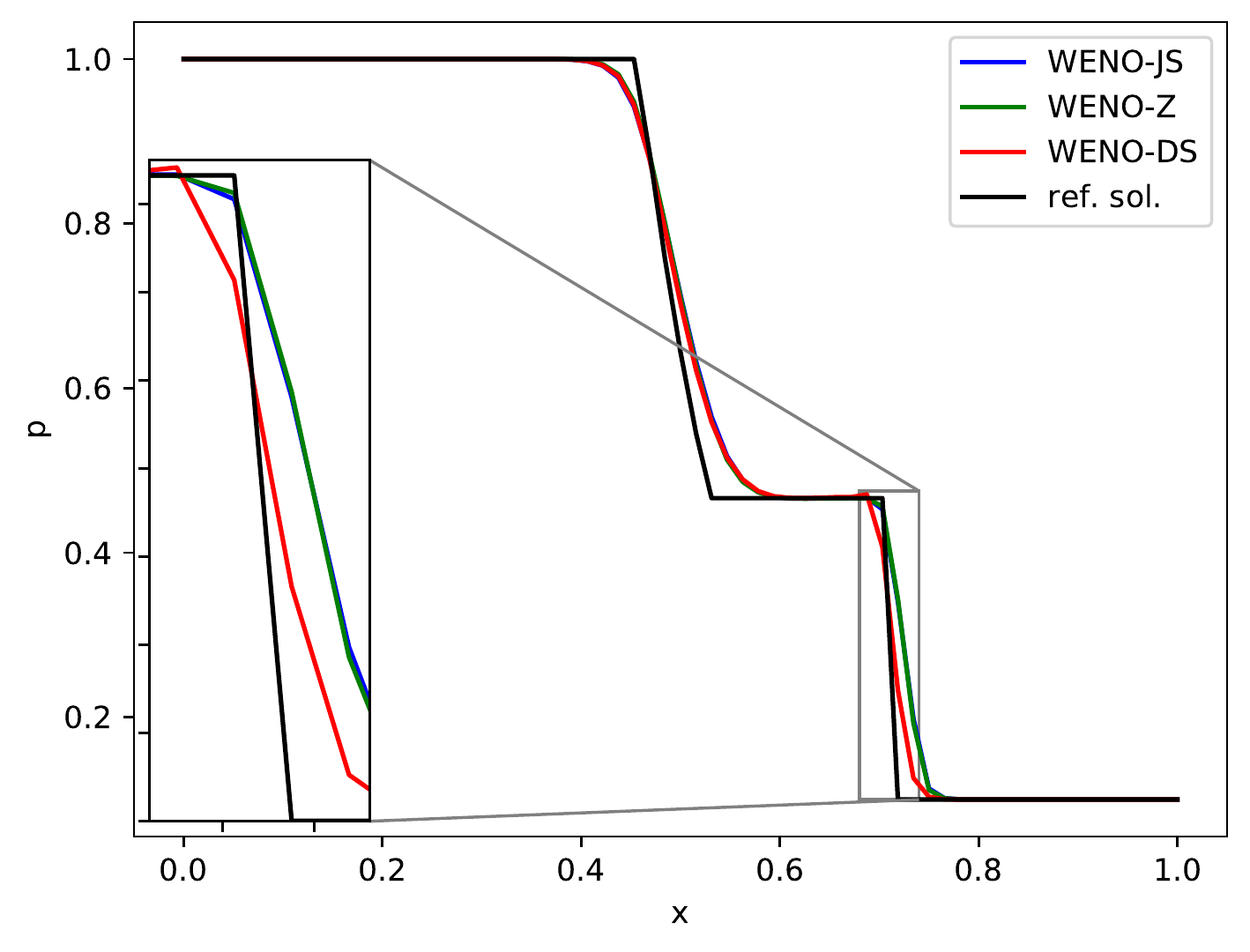} \label{fig:Sod_p_483}
\end{subfigure}
\caption{Solution of the Sod problem \eqref{eq:sod_mod}, using the WENO-JS, WENO-Z and WENO-DS obtained after the 483th training cycle, $N = 64$.} \label{fig:Sod_483}
\end{figure}

As we can see, loss decreases with the increasing number of training cycles. 
However, the final model of WENO-DS should be chosen carefully. 
In the Figure~\ref{fig:Sod_72} and \ref{fig:Sod_483} we present the solution of the modified Sod problem with the 
initial condition \eqref{eq:sod_mod} 
\begin{equation}  \label{eq:sod_mod}
  (\rho,u,p) = \begin{cases} 
    (1, 0.75, 1) \quad 0 \le x \le 0.5,\\
    (0.125, 0, 0.1) \quad 0.5 < x \le 1,\\
\end{cases} 
\end{equation}
for $\rho$, $u$ and $p$ using 64 space points and using the models obtained after the 72th training cycle and after the 483th training cycle, where the minimal value of $LOSS(\rho, u, p)$ is obtained. We compare the corresponding error values in Table \ref{tab:Sod}. 
As we can see, the second model of WENO-DS leads to much smaller errors, compared to WENO-JS and WENO-Z. 
However, we observe an oscillation on the solution for $p$ and $u$. On the other hand, 
using the first model obtained after the 72th training cycle, we see qualitatively very good solution, 
although the improvement on error is not so large.

\begin{table}[htp!]
\begin{subtable}[c]{\textwidth}
\centering
\scalebox{0.8}{\begin{tabular}{|c|c|c|c|c|c|c|c|c|}
\hline
\multicolumn{1}{|c|}{}&\multicolumn{4}{|c|}{\ $L_\infty$}&\multicolumn{4}{|c|}{\ $L_2$} \\
\hline
&  WENO-JS & WENO-Z &  WENO-DS & ratio &WENO-JS & WENO-Z &  WENO-DS & ratio\\
\hline
\hline
$\rho$  &  0.144962 &  0.145601 &  \textbf{0.117874} & 1.23&  0.022758 &  0.022120 &  \textbf{0.020395} & 1.08 \\ \hline
$p$ &  0.240541 &  0.243511 &  \textbf{0.187788}  &1.28 &  0.027821 &  0.027364 &  \textbf{0.022988} & 1.19\\ \hline
$u$ &  1.047592 &  1.055905 &  \textbf{0.886982} &1.18 &  0.107484 &  0.106745 &  \textbf{0.087345} & 1.22\\ \hline
\end{tabular}}
\caption{WENO-DS method obtained after the 72th training cycle.}
\label{tab:Sod_1}
\end{subtable}
\begin{subtable}[c]{\textwidth}
\centering
\scalebox{0.8}{\begin{tabular}{|c|c|c|c|c|c|c|c|c|}
\hline
\multicolumn{1}{|c|}{}&\multicolumn{4}{|c|}{\ $L_\infty$}&\multicolumn{4}{|c|}{\ $L_2$} \\
\hline
&  WENO-JS & WENO-Z &  WENO-DS & ratio &WENO-JS & WENO-Z &  WENO-DS & ratio\\
\hline
\hline
$\rho$ &  0.144962 &  0.145601 &  \textbf{0.105094}& 1.38 &  0.022758 &  0.022120 &  \textbf{0.018560} &1.19 \\ \hline
$p$ &  0.240541 &  0.243511 &  \textbf{0.132999} &1.81 &  0.027821 &  0.027364 &  \textbf{0.019334} & 1.42\\ \hline
$u$ & 1.047592 &  1.055905 &  \textbf{0.693131}& 1.51&  0.107484 &  0.106745 &  \textbf{0.067311} &1.59\\ \hline
\end{tabular}}
\caption{WENO-DS method obtained after the 483th training cycle.}
\label{tab:Sod_2}
\end{subtable}
\caption{Comparison of $L_\infty$ and $L_2$ error of WENO-JS, WENO-Z and WENO-DS methods for the solution of the Euler equations of gas dynamics for the Sod problem \eqref{eq:sod_mod}. As 'ratio' we denote the minimum error of the methods WENO-JS and WENO-Z divided by the error of WENO-DS (rounded to 2 decimal points).}
\label{tab:Sod}
\end{table}

\begin{figure}[htp!]
\begin{subfigure}{.32\textwidth}
\centering
\includegraphics[width=\linewidth]{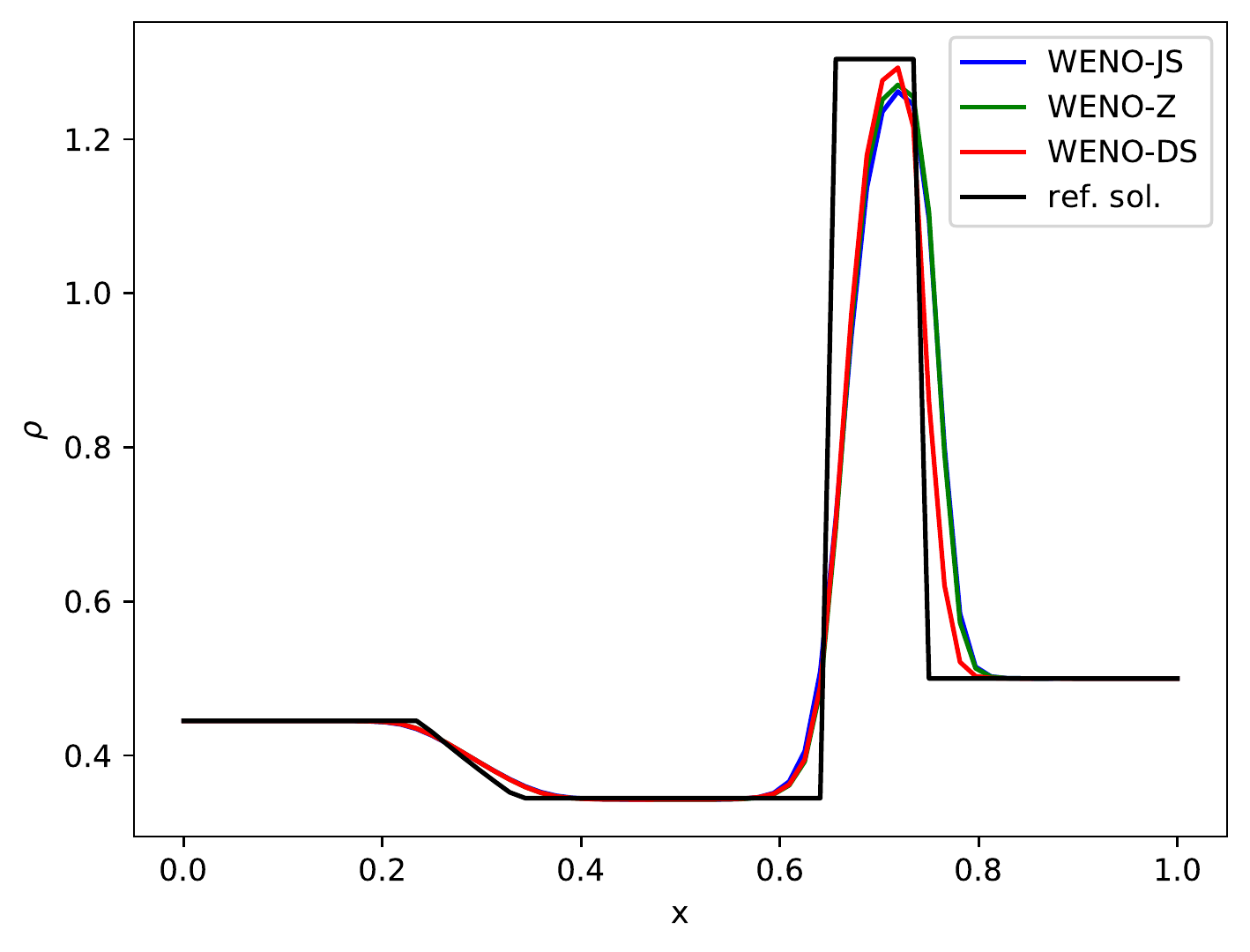} \label{fig:Lax_rho_72}
\end{subfigure}
\begin{subfigure}{.32\textwidth}
\centering
\includegraphics[width=\linewidth]{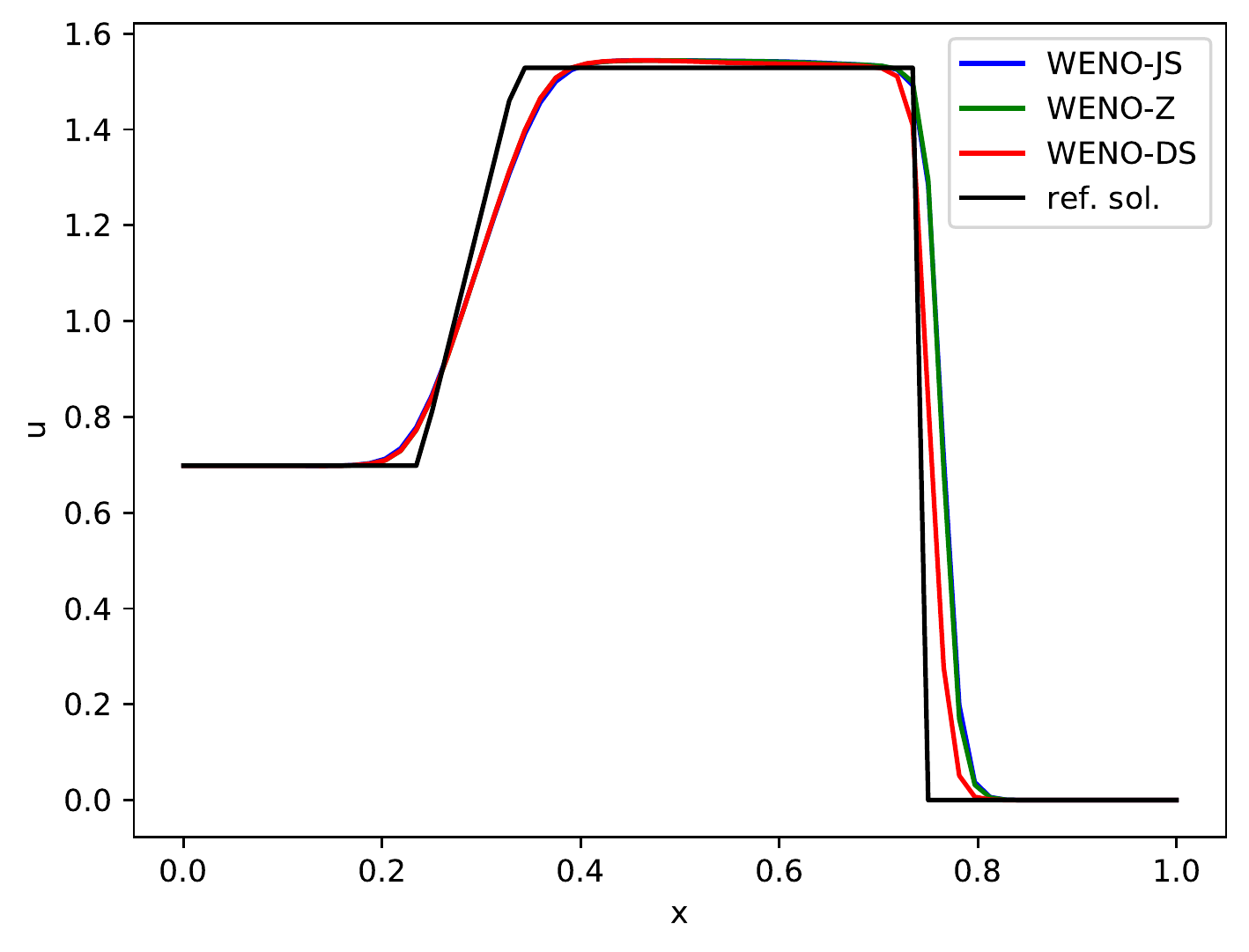} \label{fig:Lax_u_72}
\end{subfigure}
\begin{subfigure}{.32\textwidth}
\centering
\includegraphics[width=\linewidth]{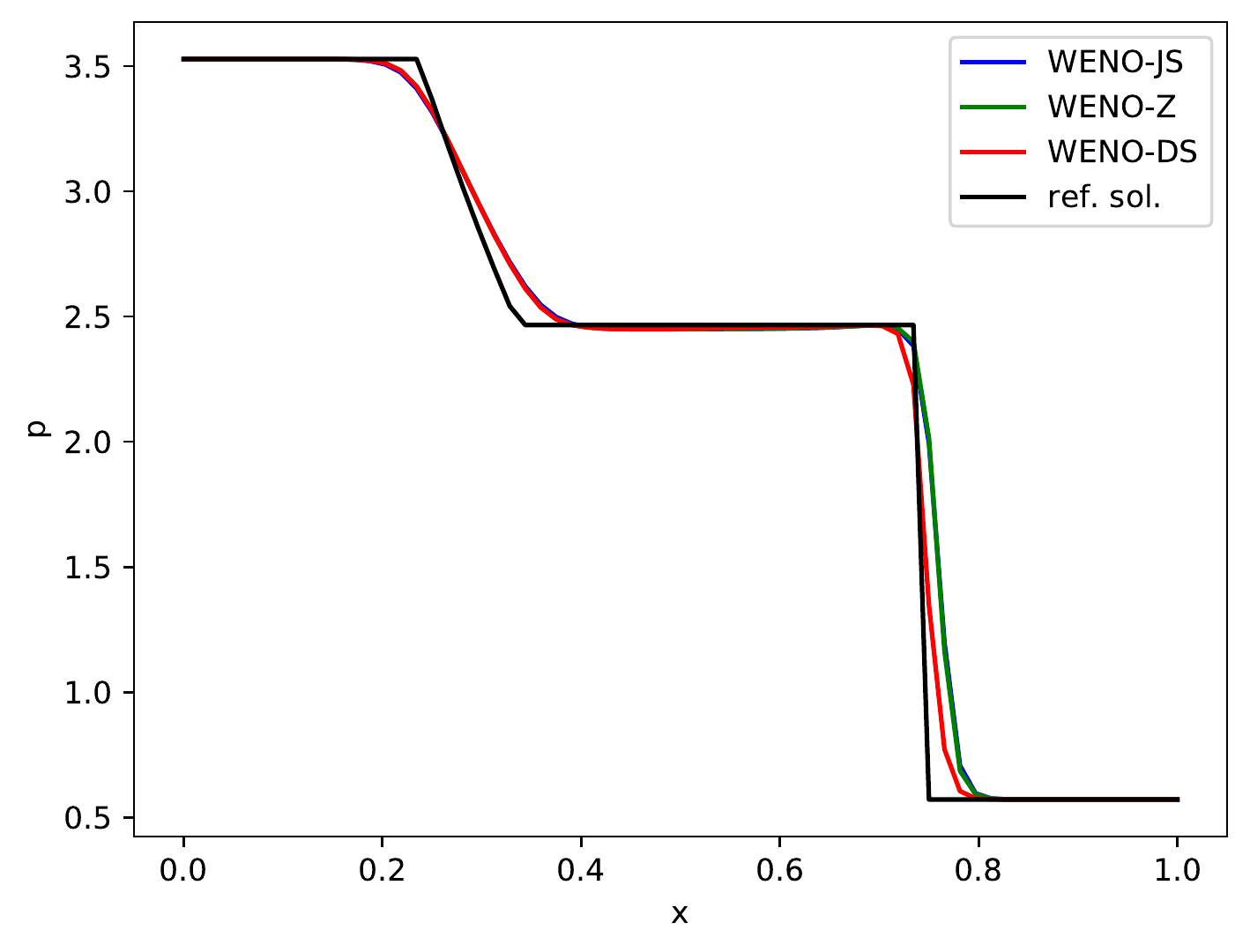} \label{fig:Lax_p_72}
\end{subfigure}
\caption{Solution of Lax problem \eqref{eq:Lax}, using the WENO-JS, WENO-Z and WENO-DS obtained after the 72th training cycle, $N = 64$.} \label{fig:Lax_72}
\end{figure}

\begin{figure}[ht!]
\begin{subfigure}{.32\textwidth}
\centering
\includegraphics[width=\linewidth]{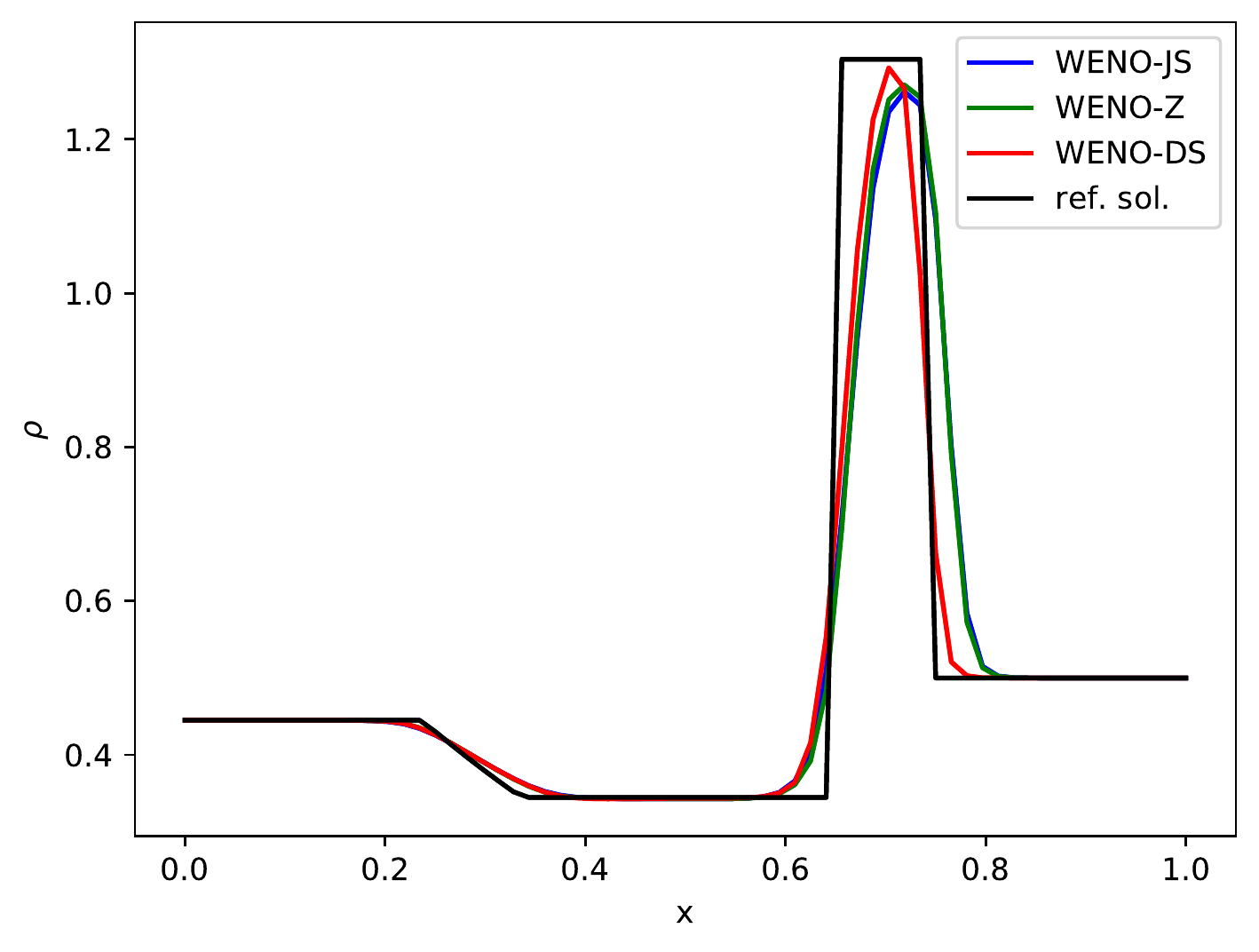} \label{fig:Lax_rho_483}
\end{subfigure}
\begin{subfigure}{.32\textwidth}
\centering
\includegraphics[width=\linewidth]{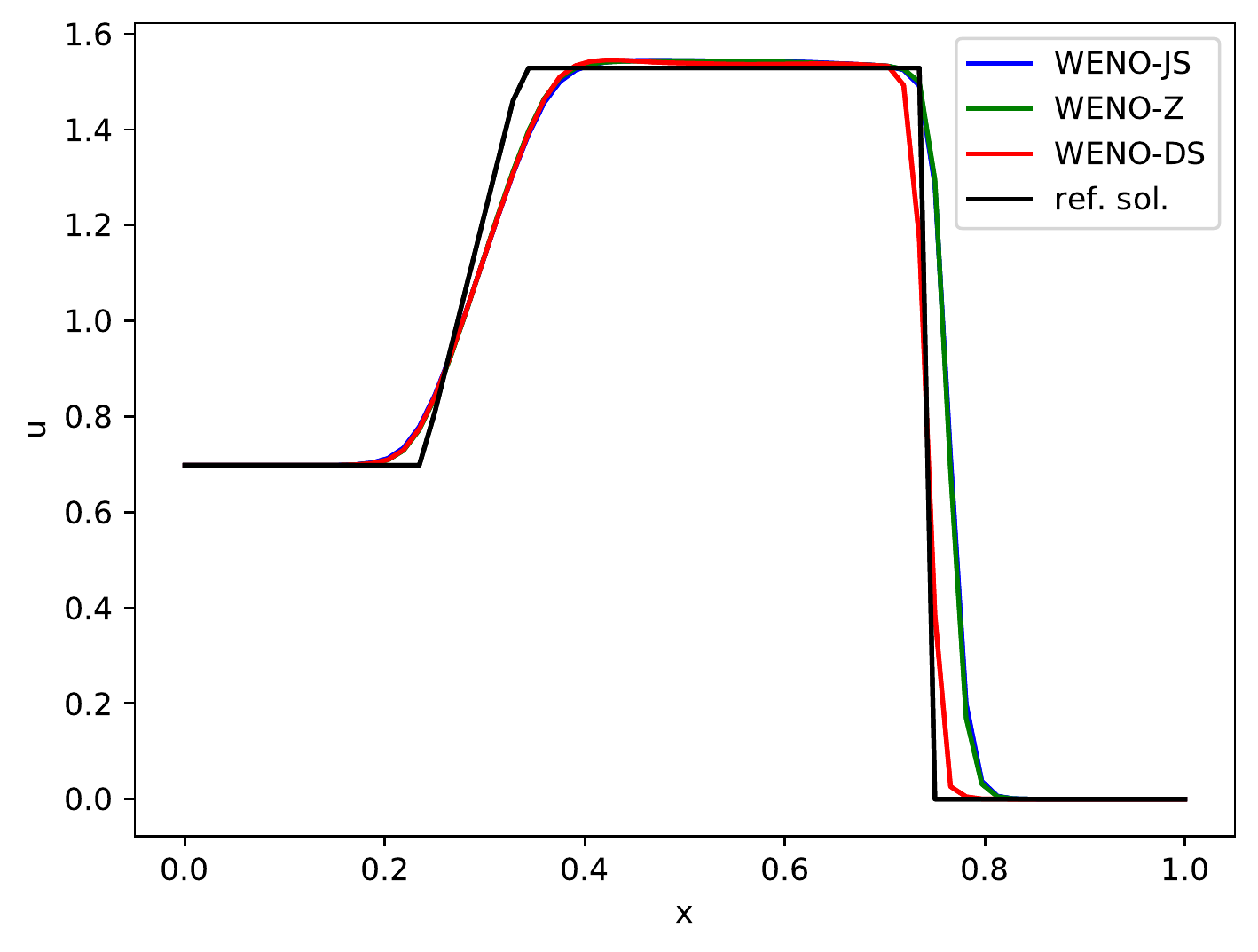} \label{fig:Lax_u_483}
\end{subfigure}
\begin{subfigure}{.32\textwidth}
\centering
\includegraphics[width=\linewidth]{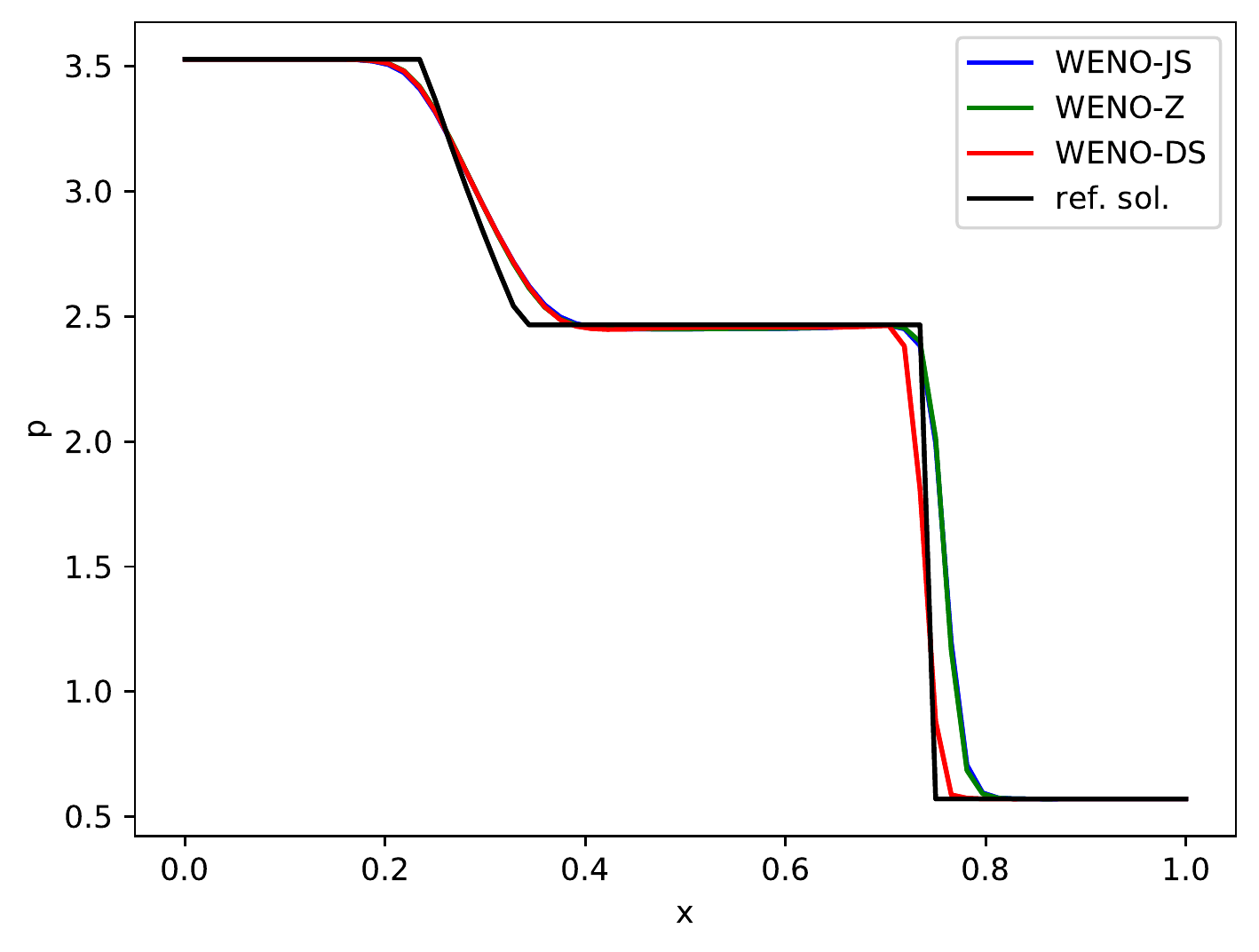} \label{fig:Lax_p_483}
\end{subfigure}
\caption{Solution of the Lax problem \eqref{eq:Lax}, using the WENO-JS, WENO-Z and WENO-DS obtained after the 483th training cycle, $N = 64$.} \label{fig:Lax_483}
\end{figure}

\begin{table}[ht!]
\begin{subtable}[c]{\textwidth}
\centering
\scalebox{0.8}{\begin{tabular}{|c|c|c|c|c|c|c|c|c|}
\hline
\multicolumn{1}{|c|}{}&\multicolumn{4}{|c|}{\ $L_\infty$}&\multicolumn{4}{|c|}{\ $L_2$} \\
\hline
&  WENO-JS & WENO-Z &  WENO-DS & ratio &WENO-JS & WENO-Z &  WENO-DS & ratio\\
\hline
\hline
$\rho$ & \textbf{0.596377} &  0.609388 &  0.599831 & 0.99&  0.088800 &  0.088287 &  \textbf{0.071987} & 1.23\\ \hline
$p$ &  1.422278 &  1.442809 &  \textbf{0.775437} & 1.83&  0.141413 &  0.141420 &  \textbf{0.079457} &1.78\\ \hline
$u$ &  1.283662 &  1.294909 &  \textbf{0.831654} &  1.54& 0.133659 &  0.132723 &  \textbf{0.081986} &1.62\\ \hline
\end{tabular}}
\caption{WENO-DS method obtained after the 72th training cycle.}
\label{tab:Lax_1}
\end{subtable}
\begin{subtable}[c]{\textwidth}
\centering
\scalebox{0.8}{\begin{tabular}{|c|c|c|c|c|c|c|c|c|}
\hline
\multicolumn{1}{|c|}{}&\multicolumn{4}{|c|}{\ $L_\infty$}&\multicolumn{4}{|c|}{\ $L_2$} \\
\hline
&  WENO-JS & WENO-Z &  WENO-DS & ratio &WENO-JS & WENO-Z &  WENO-DS & ratio\\
\hline
$\rho$ &  0.596377 &  0.609388 &  \textbf{0.507207} & 1.18 & 0.088800 &  0.088287 &  \textbf{0.061225} & 1.44\\ \hline
$p$ &  1.422278 &  1.442809 &  \textbf{0.656865} & 2.17 & 0.141413 &  0.141420 &  \textbf{0.071175} & 1.99\\ \hline
$u$ &  1.283662 &  1.294909 &  \textbf{0.386022} & 3.33 & 0.133659 &  0.132723 &  \textbf{0.052814} &2.51 \\ \hline
\end{tabular}}
\caption{WENO-DS method obtained after the 483th training cycle.}
\label{tab:Lax_2}
\end{subtable}
\caption{Comparison of $L_\infty$ and $L_2$ error of WENO-JS, WENO-Z and WENO-DS methods
for the solution of the Euler equations of gas dynamics for the Lax problem \eqref{eq:Lax}.
As 'ratio' we denote the minimum error of the methods WENO-JS and WENO-Z divided by the error of WENO-DS (rounded to 2 decimal points).}
\label{tab:Lax}
\end{table}

We also applied the trained models to the  Lax problem \cite{lax} with an initial condition
\begin{equation}  \label{eq:Lax}
  (\rho,u,p) = \begin{cases} 
    (0.445, 0.698, 3.528) \quad 0 \le x \le 0.5,\\
    (0.5, 0, 0.571) \quad 0.5 < x \le 1.\\
\end{cases} 
\end{equation}
The solution also has a left rarefaction wave, a right traveling contact wave, and a right shock wave.
In the Figure~\ref{fig:Lax_72} and \ref{fig:Lax_483} we show the solution produced with both WENO-DS models 
mentioned before and see the similar behaviour of the solution.
As indicated in Table~\ref{tab:Lax}, we observe even greater improvement on errors than in the Sod problem using both models, 
especially using the second model.

We can see, that the shape of the density profile is different, when compared to the Sod problem. 
Finally, let us note, that although we have not trained the presented models on the parameters that would lead to such a solution, the models are robust enough and can reliably detect the shocks.

\section{Conclusion}\label{concl}
In this work, we have improved the fifth-order WENO shock-capturing scheme 
by using deep learning techniques. 
To do this, we trained a relatively small neural network to obtain
modified smoothness indicators of the WENO scheme. 
This was done in a way that avoided post-processing of the coefficients to ensure consistency. 
We applied our enhancement to the WENO-Z scheme, where the fifth-order convergence on the smooth solutions can be theoretically proven.
Our new method, the WENO-DS scheme, is quite easy to use and significantly improves the numerical results, 
especially in the presence of discontinuities, even for cases that have not been trained before.
We have demonstrated our results with the inviscid Burgers' equation, 
the Buckley-Leverett equation, and the 1-D Euler equations of gas dynamics. 
Finally, let us note, that this paper can be seen as a proof of concept, that neural networks can be efficiently combined with an existing numerical scheme, preserving its convergence order.



\appendix
\section{Parameters used for generating the data set for 1-D Euler equations of gas dynamics} \label{A:1}

The problem samples representing different versions of the Euler equations of gas dynamics \eqref{eq:sod} were defined using parameters generated by the following algorithm.

\begin{algorithmic}
\STATE Choose randomly $s \in \{0,1\}$
\IF {$s = 0$} 
        \STATE $p_l = a + b$, \quad $a \in \mathcal{U}[0.5,10]$, \quad $b \in \mathcal{U}[-0.05,0.05]$, \\
        $p_r = 1/c$,  \quad $c \in \mathcal{U}[5,10]$, \\
        $\rho_l = p_l$, \\
        $\rho_r = p_r + d$,  \quad $d \in \mathcal{U}[-0.05,0.05]$, \\
        $u_l = e$,  \quad $e \in \mathcal{U}[0,1]$, \\
        $u_r = 0$,
\ELSE 
        \STATE $p_l = 1$, \\
        $p_r = 0.1$,  \\
        $\rho_l = k $,  \quad $k \in \mathcal{U}[1,3]$,  \\
        $\rho_r = \frac{1}{10} \rho_l + l$, \quad $l \in \mathcal{U}[-0.05,0.05]$, \\
        $u_l = r$,  \quad $r \in \mathcal{U}[0,1]$, \\
        $u_r = 0$,
\ENDIF    
\end{algorithmic}
where
\begin{equation} 
  (\rho,u,p) = \begin{cases} 
    (\rho_l, u_l, p_l) \quad 0 \le x \le 0.5,\\
    (\rho_r, u_r, p_r) \quad 0.5 < x \le 1.\\
\end{cases} 
\end{equation}





\bibliographystyle{elsarticle-num-names}
\bibliography{WENO.bib}







\end{document}